\input amstex
\newif\ifsyntax\syntaxtrue
\syntaxfalse
\ifsyntax\syntax\immediate\write16{**** no output: SYNTAX ONLY ****}\fi
\documentstyle{amsppt}
\magnification 1200
%       general definitions
%%% defs.tex (WJM)  standard macros for almost everything.
%	(cleaned up for sending in with SCH paper)

% 
% The definitions below are identical to those in amsppt.sty version
% 2.0 except that they are not outer.  This is necessary in order to
% use my defintion of \proc below (and in fact it is necessary for any
% way I can think of to use \proclaim with automatic numbering of theorems.
%
% I find it hard to regard the declaration of \proclaim as outer as other
% than a bug.  It was outer in early versions of version 1, but this
% was changed, presumably for this reason. I hope it will be
% officially changed back again.
%
\comment	3/17/92 (after submitted version). 
 \catcode`\@=11
%%%%%% this problem avoided by using \csname proclaim \endcsname %%%%
 \def\proclaim{\runaway@{proclaim}\def\envir@{proclaim}%
   \nofrills@{.\enspace}\proclaim@
  \DNii@##1{\penaltyandskip@{-100}\medskipamount\varindent@
    \usualspace@{{\proclaimheadfont@\enspace}}\proclaimheadfont@
    \ignorespaces##1\unskip\proclaim@
   \sl\ignorespaces}% 
  \FN@\next@}
 \catcode`\@=\active
\endcomment
\newif\ifdraft\draftfalse

%
%  standard definitions
%
\define\iwlog{{\it wlog\/}}

% greek letters (those with at most 3 letters exempted).
\let\ga=\alpha
\let\gb=\beta
\let\gd=\delta
\let\muchmore=\gg
\let\gg=\gamma
\let\gi=\iota
\let\gk=\kappa
\let\gl=\lambda
\let\gs=\sigma
\let\gS=\Sigma
\let\gth=\theta
\let\gw=\omega
\let\gz=\zeta

% script letters

\define\cc{{\Cal C}}
\define\cd{{\Cal D}}

\define\cf{{\Cal F}}

\define\ch{{\Cal H}}

%       other symbols
\define\catx#1#2{\mskip-2mu\raise#1\hbox{$#2\smallfrown$}\mskip-2mu}

	%\nothing, \restrict
\define\nothing{\varnothing}
\define\restrict{{\restriction}} %NOT spaced as a binary relation

	% orderings: \tri = triangle relation of UF's, \elem=elementary sbstr,
	%		\sat (satisfies)

\def\ps{{\Cal P}}

	% note that displaysize can't happen.  This makes the size one larger
	% than it would be otherwise.

%\define\sh{^{\tsize\sharp}}
 
\define\const#1{\frak c_{#1}}

\define\card#1{|#1|}
\define\Card#1{\Vert#1\Vert}
\define\set#1{\{\,#1\,\}}
\define\seq#1{(\,#1\,)}
%\define\ve#1{\vec{\mathstrut#1}}
\define\ve#1{\bold #1}
\define\image{\raise1.5pt\hbox{\lq\lq}\kern-.7pt} 
\predefine\integral{\int}

\define\union{\bigcup}
 
%       operators

\define\cof{\operatorname{cf}}
\define\range{\operatorname{range}}
\define\domain{\operatorname{domain}}
\redefine\sup{\operatorname{sup}}
\define\ult{\operatorname{ult}} 
\define\crit{\operatorname{crit}}
\define\len{\operatorname{len}}

%	proclaims & demos
\newcount\sectno
\global\sectno=0
\def\newsectno#1{\global\sectno=#1\global\thno=0\relax\the\sectno\ ---}
\outer\def\secthead#1#2{
	\global\thno=0\global\sectno=#1\subheading{\the\sectno.\ #2}}

%\outer\def\secthead#1{
%        \advance\sectno by 1
%        \subheading{\number\sectno.\ #1}
%}
\newcount\thno
\global\thno=0
\def\thethm{\the\sectno.\the\thno}
\def\proc#1{
	\advance\thno by 1
	\csname proclaim\endcsname{\thethm\ #1}
}
\outer\def\cor{\proc{Corollary}}
\outer\def\theorem{\proc{Theorem}}
\outer\def\lemma{\proc{Lemma}}
\outer\def\prop{\proc{Proposition}}
\outer\def\defin{\proc{Definition}}
\define\defining#1{{\it #1}}

\let\endth=\endproclaim

	%% thmtags
% \writetagfile{<name of file>}
	% -> also write "thmfile.tex", must follow \usethmtags
% \thmtag\<word> -> make \<word> = current \sectno.\thno
% use \input <file> before \writetagfile

\def\putinmargin#1{\leavevmode\vadjust{
   \vbox to0pt{\vss\hbox to 0pt{\hss\hbox to0.70truein{\eightpoint#1\hss}}}}}
\newwrite\thmfile

%% default definition \thmtag if \writetagfile not used:
\def\thmtag#1{\edef#1{\number\sectno.\number\thno}
	\ifdraft{\escapechar=-1\putinmargin{\string#1}}\fi}

\def\nextsect#1/{{\def\xxx#1{\count0 by #1}\count0=\sectno
	\expandafter\advance\xxx#1\the\count0}}
\def\nextthm#1/{{\def\xxx#1{\count0 by #1}\count0=\thno
      \expandafter\advance\xxx#1\the\sectno.\the\count0}} %use backthm instead.
\def\backthm#1{{\def\xxx#1{\count0 by -#1}\count0=\thno
	\expandafter\advance\xxx#1\relax\the\sectno.\the\count0}}
\define\proof{\demo{Proof}}
\define\endproofof#1{%
  \nobreak\hfil\penalty0\hfilneg\quad\qed\nobreak\quad\nobreak(#1)%
  \enddemo}
%\define\endproofof#1{\newline\qed\qquad(#1)\enddemo}
\define\endproof{\qed\enddemo}

\expandafter\ifx\csname today\endcsname\relax
	\def\today{\number\month/\number\day/\number\year}\fi
 
\def\lineto#1/{
	\def\xxx#1{height .4pt depth 0pt width #1}
	\expandafter\vrule\xxx#1
}

%----------------------------------------------------------------
% marginal notes
\define\squiggle{{\def\m{\muchmore}$\m\m\m\m\m$}}
\define\markline{{\leavevmode\vadjust{%
	\vbox to 0pt{\vss\hbox to 0pt{\hss\raise.25ex\hbox{\squiggle\ \ }}}}}}

%\preliminary\immediate\write16{**** PRELIMINARY ****}
\def\mainthm{1.1}
\def\CoveringLemma{1.2}
\def\hsupxiLem{2.3}
\def\wksupportLem{2.4}
\def\indiscDef{2.5}
\def\indiscLem{2.6}
\def\cXeqcXpr{2.7}
\def\semicont{2.8}
\def\accptDef{2.9}
\def\supportOne{2.10}
\def\havesupports{2.11}
\def\unique{2.13}
\def\uniqueFinally{2.14}
\def\mainthree{3.1}
\def\hCclosureDef{4.2}
\def\haveStdSys{4.3}
\def\XisMax{4.4}

\define\core{K}
\define\coref{\core(\cf)}
\define\barcore{\overline K}
\def\ch{\Cal H}

\define\vc{\ve c}
\define\vd{\ve d}
\define\barcf{\overline{\Cal F}}
\define\barcc{\overline{\Cal C}}
\define\barm{\overline{m}}
% \define\claim{\proclaim{Claim}}    %%% fixed 3/17/92
   \define\claim{\csname proclaim\endcsname{Claim}}
\redefine\ve#1{\vec#1}
\define\ccx{\cc^X}
\define\cohere{\frak C}

\define\fnote#1{\plainfootnote{{\bf/*/}}{#1}}
%\define\changed{\plainfootnote{{\bf/**/}}{changed!}}
\define\changed{\write16{Changes on page \the\pageno.}}
%________________________________________________________________ 
\topmatter
\title On the Singular Cardinal Hypothesis\endtitle
\author W. J. Mitchell\endauthor
\thanks{
Some of the work in this paper was done while the author was visiting
the Hebrew university with support from the Lady Davis Foundation, and
while the author was visiting UCLA and the
California Institute of Technology.
\endgraf\noindent
This work was partially supported by grant number DMS-8614447 from the
National Science Found\-a\-tion.}
\endthanks
\date{\today\ifdraft -- preliminary\fi}\enddate

 \address{Department of Mathematics, University of Florida,
 Gainesville, Florida, 32611}\endaddress
 \email mitchell\@math.ufl.edu\endemail

\endtopmatter
%________________________________________________________________ 
\document
%% \openup2\jot	%% for referee, remove for printing

% \ifsyntax\else\writetagfile{sch-tags.tex}\fi
	% comment if texing and printing out only part of the paper.
\heading{\newsectno1 Introduction}\endheading

The Singular Cardinal Hypothesis (SCH) asserts that if $\gk$ is any singular
strong limit cardinal then $2^{\gk } = \gk ^{+} $.  It is known to be
consistent that the SCH fails: Prikry \cite{Pr} obtains a model of $\lnot$SCH
from a model in which the GCH fails at a measurable cardinal $\gk $, and
Silver in turn (see \cite{KM}) obtains the failure of the GCH at a measurable
cardinal from a model for a $\gk ^{+} $-supercompact cardinal. 
Silver's result has been improved by Woodin \cite{W} and again by Gitik
\cite{G} so that the consistency of ``$\gk$ measurable and $2^{\gk } >\gk ^{+}
$'', and hence of $\lnot $SCH, is now known to follow from that of $o(\gk
)=\gk ^{++} $.  Results in \cite{Mi?} show that the failure of the GCH at a
measurable cardinal is in fact equiconsistent with $\exists \gk\,o(\gk )=\gk
^{++} $. 

The best previously published lower bounds on the consistency strength of
$\lnot $SCH are those of Dodd and Jensen, \cite{D} which use the covering
lemma for $L^{\mu }$ to show that at least one measurable cardinal (and
slightly more) is required.  In this paper we use the core model $K(\Cal F)$
for sequences of measures (\cite{Mi84b}, \cite{Mi?}, \cite{Mi87}) to prove:

\theorem\thmtag\mainthm (i) Con($\lnot $SCH) implies Con($\exists \gk \forall
\ga <\gk \exists \nu <\gk\;o(\nu )\ge \ga $).

  (ii) If $\gk$ is a singular strong limit cardinal with $\cof(\gk )>\gw$ and
$2^{\gk } >\gk ^{+}$ then there is an inner model with $o(\gk )= \gk ^{++} $. 
\endth

Since this paper was originally
submitted Gitik has improved this result to give the exact
consistency strength of the failure of the SCH. In 
\cite{G?} he strengthens the conclusion of
clause~(i) to $o(\gk)=\gk^{++}$, which is shown in
\cite{G} to be best possible. The arguments of \cite{G?} use results
from  section~2 of this paper but his proof, which uses results of
Shelah, is substantially different from our arguments in section~4.
Our proof gives some additional information in the limited situation
in which it applies and may be of some independent interest.
In \cite{G??} he strengthens the conclusion
of theorem~\mainthm(ii) to $o(\gk )= \gk^{++} +\gl $, where $\gl =
\cof(\gk)$.  This result is best possible, as is
shown by a generalization by McDermitt \cite{Mc} of Woodin's work. In
this case Gitik does use  our argument from section~3 of this paper,
with one new technique which overcomes the difficulty encountered in
working with extenders in place of measures.

Our proof of theorem~\mainthm\ depends on a refined version of the
covering lemma for sequences of measures from \cite{Mi?, section~6}.
Theorem~\CoveringLemma\ below is the  basic statement of this
refinement, but we will need more detailed information about the
system $\Cal C$ of indiscernibles.  This information appears in
section~2 along with the proof of theorem~\CoveringLemma. The
covering set $h\image (\gg ;\Cal C)$ in the statement of the lemma is
the  smallest set containing $\gg $, closed under $h$, and containing
$\Cal C(\ga,\gb)$ whenever it contains $\ga $ and $\gb $. 

\theorem\thmtag\CoveringLemma Assume that there is no inner model of $\exists \gk\
o(\gk )=\gk^{++}$. Then for any set $w$ of ordinals there is a function 
$h\in K(\Cal F)$, an ordinal $\gg <\card{w^{\gw }}^{+} $, and a system 
$\cc$ of 
indiscernibles for $K(\Cal F)$ such that $h\image(\gg ;\Cal C)\supset w$. 
\endth

This version is much closer to the covering lemmas for $L$ and $L^{\mu
} $, proved by Jensen and by Jensen and Dodd respectively, than that
given in \cite{Mi?}, but it is weaker than those versions in one
crucial respect: a direct generalization of their results would give a
single system $\Cal C$ which works for all sets $w$.  This
generalization is known to be false, (see the discussion following the
proof of theorem~4.1 in section~4) but the system $\Cal C$ given by
theorem~\CoveringLemma\ is, in a sense to be made precise later, unique and
maximal up to finite changes for the measures which lie in the
covering set $h\image(\gg ;\Cal C)$.  The proof of theorem~\mainthm\
will depend on using this maximality and uniqueness to define a small
set of ``standard'' systems of indiscernibles which is large enough
that theorem~\CoveringLemma\ is still true if the system $\Cal C$ is required to
be taken from this set. 

The Singular Cardinal Hypothesis can be stated in a more general form
which does not assume that $\gk$ is a strong limit cardinal: If $\gk$
is any singular cardinal then
$\gk^{\cof(\gk)}=\max\left(\gk^{+},2^{\cof(\gk)}\right)$.  Dodd and Jensen's
results in \cite{D} are still valid for this stronger statement of the
SCH.  The same is probably true of theorem~\mainthm\ but a proof would
require using techniques of Dodd and Jensen to avoid the assumption,
frequently used both in \cite{Mi?} and in this paper, that the set
$X\prec H_{\gk^{++}}$ used  in the proof of the covering lemma is
closed under $\gw$ sequences.

The  proof  of  theorem~\CoveringLemma,  together  with  the  analysis  of  the
maximality  of  the  sequence  of  indiscernibles,  is  in  section~2.
Section~3  contains  the  proof  of  theorem~\mainthm(ii),  the  case
of uncountable  cofinality,  and  section~4  contains  the  proof  of
theorem~\mainthm(i),  the  case  of  countable cofinality.
Each of  the  sections~3
and~4  depends on  section~2,  but  they  can  be  read  independently  of  each
other.

This  paper  depends  heavily  on  the  results  and  methods  of
\cite{M84b}  and  \cite{M?}.  We  have  attempted  to  summarize  the  necessary
facts  in  this  paper,  but  an  acquaintance  with  that  material  would  be
desirable.  At  the  least  an  understanding  of  the  theory  of  coherent
sequences  of  measures  (\cite{Mi74}, \cite{Mi83})  is  necessary.  We  summarize
below  some  of  the  notation  from  these  sources.

\subheading{Notation} A {\it sequence of measures} is a coherent 
function $\Cal F$ with domain of the form
$\set{(\ga ,\gb ):\ga <l^{\Cal F}\text{ and } \gb <o^{\Cal F} (\ga )}$, 
where $l^{\Cal F} $ is an ordinal and $o^{\Cal F} $ is a function,
such that $\cf(\ga,\gb)$ is a measure on $\ga$ for all ordinals 
$\ga<l^\cf$ and $\gb<o^\cf(\ga)$.
The requirement that $\cf$ be coherent means that for all ordinals $\ga$,
$\gb'$ and $\gb$ such that
$\ga <l^{\Cal F} $ and $\gb '<\gb <o^{\cf}(\ga )$, if $f\in\coref$ is the
least function such that $[f]_{\cf(\ga,\gb)}=\gb'$ then 
for all $x\in\Cal P(\ga)\cap\coref$ we have $x\in
\Cal F(\ga ,\gb ')$ iff $\set{\nu :x\cap \nu \in \Cal F(\nu ,f(\nu ))}\in \Cal
F(\ga ,\gb )$.
We write $\cohere(\ga ,\gb ',\gb )$ for this {\it coherence function} $f$. 
The coherence of $\cf$ implies that if
$i:\coref\to\ult(\coref,\cf(\ga,\gb))$ is the canonical embedding
then $i(\cf)(\ga,\gb')=\cf(\ga,\gb')$ for all $\gb'<\gb<o^{\cf}(\ga)$.

The precise definition of a system of indiscernibles is given in the
next section as
definition~\indiscDef, but we present here a review of the basic
theory. A  {\it
system of indiscernibles} for $\Cal F$ is a function $\Cal C$ with
$\domain(\cc)\subset\domain(\cf)$ such that  $\Cal C(\ga,\gb)$ is a
subset of $\ga$ whenever $(\ga,\gb)\in\domain(\cc)$. As one would
expect, the sets $\cc(\ga,\gb)$ resemble Prikry sequences in the sense
that if $\cc(\ga,\gb)$ is cofinal in $\ga$ then  $$\forall x\in\Cal
P(\ga)\cap\coref\;\bigl( x\in\Cal F(\ga ,\gb )\iff
(\cc(\ga,\gb)\setminus x\text{ is bounded in $\ga$})\bigr)$$
but the ordinals in $\cc$ also work uniformly as indiscernibles, even
when they belong to different measures:
if $\vec
c$, $\vec\ga$, and $\vec\gb$ are $\gw$-sequences such that $\vec c$ is
increasing,  $c_i\in\cc(\ga_i,\gb_i)$ for all $i\in\gw$, and
$\union_i c_i=\union_i \ga_i$, then for any function $g\in\coref$
there is $i_0<\gb$ such that for all $i>i_0$ $$\forall x\in g\image
c_i\;(c_i\in x\iff x\in\cf(\ga_i,\gb_i)).$$

The systems of indiscernibles used in this paper will come, directly
or indirectly, from
iterated ultrapowers.
Since we are making one minor change
from \cite{Mi84b} in the construction of a system of
indiscernibles from an iterated ultrapower, we give the full definition
here. Suppose, in general, that 
 $i_{0,\gth}:M_0\to M_{\gth}$ is an iterated ultrapower, with
$M_{\ga+1}=\ult(M_\ga,\cf_\ga(\gk_\ga,\gb_\ga))$ for each $\ga<\gth$
where $\cf_0$ is a coherent sequence of measures in $M_0$,
$\cf_\ga=i_{0,\ga}(\cf_0)$, and the sequence $\gk_\ga$ is strictly
increasing.
For $\gk\in \range(i_{0,\ga})$, we
define
$$\cc_\ga(\gk,\gb)=\set{\gk_{\ga'}:\gk=i_{\ga',\ga}(\gk_{\ga'})\text{ and
}\gb=i_{\ga',\ga}(\gb_{\ga'})}.$$
If $\gk\notin\range i_{0,\ga}$ then $\cc_\ga(\gk,\gb)$ is defined to
be empty unless $\gk\in\cc_\ga(\gl,\gg)$ for some
$\gl\in\range(i_{0,\ga})$, 
in which case let $\gb'=i_{\ga',\ga}(\gb)$ where $\gk=\gk_{\ga'}$,
so that
$\gb'<\gg$ and  $\gb=\cohere(\gl,\gb',\gg)(\gk)$. Then we set
$$\cc_\ga(\gk,\gb)=\cc_\ga(\gl,\gb')\cap(\gk\setminus\gk_\xi)$$ where
$\xi<\ga$ is the least ordinal such that $\gg\in \range(i_{\xi,\ga})$.
The definition given in \cite{Mi84b} was equivalent except that $\gk_\xi$
was omitted, so that in the second case $\cc_\ga(\gk,\gb)$ would be
$\cc_\ga(\gl,\gb')\cap\gk$.

Suppose that $\cc=\cc_{\gth}$ is a system of indiscernibles constructed as above.
Some basic consequences of this construction are that
$\union_{\gb}\cc(\ga,\gb)$ is closed in $\ga$ and that if
$\gb\not=\gb'$ then $\cc(\ga,\gb)\cap\cc(\ga,\gb')=\nothing$.
If $c$ is an indiscernible in $\cc$ then there will be a unique pair
$(\ga,\gb)$ such that $\ga$ is not an indiscernible and
$c\in\cc(\ga,\gb)$.  We will write $\ga^\cc(c)$ and $\gb^\cc(c)$ for this
pair of ordinals.  These ordinals have another equivalent
characterization: if $c=\gk_\nu$ in the iteration then $\ga^\cc(c)$ is
the least member of $\range(i_{\nu,\gth})\setminus c$ and $\gb^\cc(c)$ is
the unique member $\gb$ of $\range(i_{\nu,\gth})\cap o(\ga)$ such
that for all $x\in\range(i_{\nu,\gth})\cap\ps(\ga^\cc(c))$ we have $c\in
x\iff x\in\cf_{\gth}(\ga^{\cc}(c),\gb)$.

Now suppose that $m_0$ is a mouse with projectum $\rho$ and
$i=i_{0,\gth}:m_0\to m=m_\gth$ is an iterated ultrapower such that
$i\restrict\rho$ is the identity. Let $h_0$ be the canonical skolem
function for $m_0$. Then $i$ maps $h_0$ to the canonical skolem function
$h$ for $m$, so that $h\image\rho=\range(i)$ and in general if
$c=\gk_\nu$ then $h\image c=\range(i_{\nu,\gth})$. Hence the
characterization of the functions $\ga^\cc$ and $\gb^\cc$ may be
given in a form that depends, apart from their domain, only on
the function $h$; namely $\ga^\cc(c)=\min(h\image c\setminus c)$ and
$\gb^\cc(c)$
is the unique ordinal $\gb\in h\image c\cap o(\ga^\cc(c))$ such that 
$$\forall x\in h\image c\cap\ps(\ga^\cc(c))\;\bigl(c\in x\iff
x\in\cf(\ga^\cc(c),\gb)\bigr).\tag1$$

Following this idea, we use this last characterization of $\ga^\cc$
and $\gb^\cc$ to define functions $\ga^m(\nu)$ and $\gb^m(\nu)$  for
all ordinals $\nu$ in $m$ rather than only for indiscernibles $c$ in
$\cc$. Thus the functions $\ga^m$ and $\gb^m$ are definable over the
mouse $m$, but the functions $\ga^\cc$ and $\gb^\cc$ are equal to
$\ga^m\restrict C$ and $\gb^m\restrict C$, respectively, where
$C=\union_{\ga,\gb}\cc(\ga,\gb)$.  In the special case when
$\ga^m(\nu)=\nu$ or there is no ordinal $\gb$ satisfying formula~(1)
we set $\gb^m(\nu)=o^m(\nu)$. In particular $\ga^m(\nu)=\nu$ and
$\gb^m(\nu)=o(\nu)$ whenever $\nu\in h\image\nu$. 

The system of indiscernibles $\cc$ given by the covering
lemma~\CoveringLemma\ will be constructed indirectly from an iterated
ultrapower, and  as a result the idea of the last paragraph will apply,
except that $m$ will be a $\cf\restrict\gk$-mouse and hence a member
of $\coref$.   The construction
of a system of indiscernibles from an iterated ultrapower was changed
from that of \cite{Mi84b} so that
$\cc$ would be a system of indiscernibles for $\coref$ as well
as for $m$. To see why the change is necessary, note that since $m\in\coref$ the functions $\ga^m$ and $\gb^m$ are
members of $\coref$. Now define a function  $k$ on the ordinals of $m$
by letting $k(\nu)$ be the least ordinal $\gd<\nu$ such that
$\gb^m(\nu)\in h\image\gd$. Then $k\in\coref$, and hence
indiscernibility requires that (with at most finitely many exceptions)
we have $k(\ga)<c<\ga$ whenever $c\in\cc(\ga,\gb)$.
The change in the construction insures that this
requirement is satisfied.

If $x$ is any member of
$m$ then $x$ has a minimal support $d$ in the iterated ultrapower, that is,
there is a finite set $d$ of indiscernibles from the system
$\cc$ generated by $i$ such that $x\in h\image(\rho\cup d)$, and such
that $d$ is contained in any other set $d'$ of indiscernibles such that 
$x\in h\image(\rho\cup d')$.

The expression $K(\Cal F)$ will always denote the maximal core model
for sequences of measures as defined in \cite{Mi?}.  A
$\cf\restrict\gk$-mouse is a model $m=J_\xi(\cf^m)$, with
$\cf^m\restrict \gk = \cf\restrict \gk $ and $o^\cf(\ga)<\gk$ for $\ga<\gk$, 
such that $m$ is iterable
and every member of $m$ is definable in $m$ from parameters in $\gk
\cup p^m$ for some finite set $p^m$ of ordinals.  The ordinal $\gk $
is refered to as the {\it projectum} of $m$.  It should be noted that the
coherence function $\cohere^m$ for $\cf^m$ will be different from the
coherence function $\cohere$ for $\cf$, but since
$\cf^m\restrict\gk=\cf\restrict\gk$ we have
 $\cohere^{m} (\ga ,\gb ',\gb )=
\cohere(\ga ,\gb ',\gb )$ whenever $\ga<\gk$.
We will not be using $\cohere^m(\ga,\gb',\gb)$ for ordinals
$\ga >\gk $, so it follows that the problem only arises for $\ga=\gk$.

For the benefit of those with some acquaintance with the fine structure
sections of \cite{Mi84b} it should be acknowledged that some of our
discussion has been somewhat sloppy. In the notation from that paper,
a mouse $m_0$ is actually the $\Sigma_{n}^{*}$-code of a structure of the form
$J_{\xi_0}(\cf^{m_0})$, for some $n\in\gw$, and an iterated ultrapower
of $m_0$ is the internal ultrapower of the $\gS^*_n$-code $m_0$. This iterated
ultrapower may be regarded as a $\gS^*_n$-ultrapower of
$J_{\xi_0}(\cf^{m_0})$ which has all of its critical points 
$\gk_\ga$ in the interval between the $n+1$st and $n$th
projectums of $J_{\xi_\ga}(\cf^{m_\ga})$.

We write $H_{\gk}$ for the sets hereditarily of cardinality less than $\gk$,
and if $M$ is a model of set theory then we write $H^{M}_{\gk}$ for $H_{\gk}$
as defined in $M$.

\heading{\newsectno2 The Covering Lemma}\endheading

In this section we will prove theorem~\CoveringLemma, together with
various results giving more information about the system $\cc$ of
indiscernibles.  For background, and in order to provide a framework
and notation for our work, we begin with an outline of the proof of
the covering lemma from \cite{Mi?}.  Many of the ideas of this section
were previously used in \cite{Mi87}.
The results from this section will be used in the rest of
the paper to show that there is a collection of at most $\gk^+$ 
systems of indiscernibles which is rich enough to cover all
small subsets of $\gk $.
The Singular Cardinal Hypothesis will then follow from the fact that
most $\gk^{+}$ sets can be covered using only $\gk^{+}$ systems.

Let $w$ be the set to be covered, let $\gk=\sup(w)$, and assume
\iwlog\ that $\vert
w\vert^{\gw }<\vert \gk \vert $.  We can also assume \iwlog\ that $\gk $ is a
cardinal in $\coref$.  Instead of working with $w$ directly, pick a set 
 $X\supset w$ such that $X\prec H_{\gk^{++}}$, ${^\gw}X\subset X$,
and $\card X=\card{w}^{\gw}$. 
We will obtain a covering of $X$ which satisfies all of the conditions
of
theorem~\CoveringLemma\ except that the condition $\gg<\card{w^\gw}^+$
is replaced by the weaker condition
$\gg<\gk$. Theorem~\CoveringLemma\ follows from this weaker result by
a simple induction.

The first part of this section outlines the construction from
\cite{Mi?} of the $\cf\restrict\gk$-mouse $m^X$ and system
$\cc^{X}$ of indiscernibles for $m$ such that $X$, and hence $w$, can be
covered using $m$ and $\cc^{X}$. The rest of the section will contain
results showing that $\cc^X$ is a system of indiscernibles
for $\coref$ rather than just for the mouse $m$, and finally that
the system $\cc^X$ is, in an appropriate sense,
unique and maximal.

\medskip
The basic proof of the covering lemma, as given in \cite{Mi?}, is as
follows:  Let $N$ be the transitive collapse of $X$ and let
$\pi:N\cong X\prec H_{\gk^{++}}$ be the isomorphism. We will use an
overbar to indicate a preimage under the map $\pi$, so that for
example $\pi(\bar\gk)=\gk$ and $\pi(\barcf)=\cf\restrict\gk$.  Since
$\card X<\card\gk$ and $X$ is cofinal in $\gk$, $\pi$ is not the
identity on $\bar\gk$.  Now let $\gd$ be the critical point of $\pi$.
If $\ps(\gd)\cap\coref\subset N$ then
$U=\set{x\in\ps(\gd)\cap\coref:\gd\in\pi(x)}$ is an ultrafilter which
is not in $\coref$.  Now $U$ is countably complete since
${^\gw}N\subset N$, so $U$ could have been included in the sequence
$\cf$. The only reason why it would not have been included is if
$\coref$ already satisfied that $o^{\cf}(\gd)=\gd^{++}$, which would
contradict the assumption that there is no model with such a cardinal
$\gd$.
It follows that $\ps(\gd)\cap\coref\not\subset N$ and since
$\cf\restrict\gd=\barcf\restrict\gd$ it follows that there is an
$\barcf\restrict\gd$-mouse which is not in $N$.  Let $\barm_0$ be the
least mouse not in $N$, that is, $\barm_0$ is a
$\barcf\restrict\xi$-mouse for some $\xi<\bar\gk$ and there is no
smaller $n$ which is a $\barcf\restrict\xi'$-mouse for any
$\xi'<\gk$. Then there is an iterated ultrapower
$i_{0,\theta}:\barm_0@>>>\barm_{\theta}=\barm$ so that
$\cf^{\barm}\restrict\bar\gk=\barcf$.  The minimality of $\barm_0$ implies
that $(H_{\bar\gk})^{\barm}\subset N$, and the embedding $i_{0,\gth}$
generates a system $\barcc$ of indiscernibles for $\barm$.

Because $(H_{\bar\gk})^{\barm}\subset N$ we can extend
$\pi\restrict\bar\gk$ to a map ${\pi^{*}:\barm @>>> m}$ where
$m=m^{X}$ is an $\cf\restrict \gk$-mouse.   Again using
${^\gw}N\subset N$, $m$ is iterable and hence is a member of $\coref$.
Then the system $\cc^{X}$ defined by setting
$\cc^{X}(\pi^{*}(\ga),\pi^{*}(\gb))=\pi\image\big(\barcc(\ga,\gb)\big)$
is a system of indiscernibles for $m^{X}$, and the canonical skolem
function $\bar h$ for $\barm$ maps to the canonical skolem function
$h^{m}$ for $m^{X}$.  Let $\rho=\pi(\bar\rho)$.  Then $\rho<\gk$. Also
$\barm=\bar h\image(\bar\rho;\barcc)$ since $\bar\rho$ is the
projectum of $\barm$, and it follows that $w\subset
X\cap\gk=\pi\image(\barm\cap\bar\gk)=h^{m}\image\left(\pi\image
(\bar\rho);\cc^X\right)\cap\gk\subset
h^{m}\image\left(\rho;\cc^X\right)$.

If $d\subset\gk$ then we will generally abuse notation by writing
$h^m\image(d)$ for  $h\image[\rho\cup d]^{<\gw}$.

The projectum of $m^{X}$ is equal to $\gk$ rather than
$\rho$, and in fact $H_{\gk}\cap\coref\subset m$.
To see this note that $\Cal P(\bar\gk)\cap
N\cap\core(\barcf)\subset\barm$
since $\barm$ is larger than any mouse which is a member of $N$. But
$\pi^{*}$ is equal to $\pi$ on $\Cal P(\bar\gk)\cap N$, so 
$\ps(\gk)\cap\coref\cap X\subset m$.  
The elementarity of $X$ implies that it contains a subset $x$ of
$\gk$ which is is constructed later in the canonical order of construction of
$\coref$ than any bounded
subset of $\gk$.  It follows that $x\in m$ and hence every bounded subset of
$\gk$ in $\coref$ is in $m$.

\medskip
One complication should be mentioned here.  The system $\cc^{X}$ is a
system of indiscernibles for the sequence $\cf^{m}$ of measures in
$m$, rather than for the sequence $\cf$ of $\coref$.  This means that
the domain of $\cc^X$ is equal to the domain of $\cf^m$, rather than
to the domain of $\cf$. This doesn't matter below $\gk$, since
$\cf^{m}\restrict \gk=\cf\restrict\gk$, but if $o(\gk)\ge\gk^+$ then
$\cf^m$ and $\cf$ will differ at $\gk$. Let $A$ be defined in $\coref$
as  $\set{\nu:o(\nu)=\nu^{+}}$.   Then  the unique $\eta$ such that
$A\in\cf^m(\gk,\eta)$ will be  ${\gk^+}^m$ and the unique $\eta$ such
that $A\in\cf(\gk,\eta)$ will be  ${\gk^+}$.  Thus
$\cf^m(\gk,{\gk^+}^m)$ extends naturally to   $\cf(\gk,\gk^+)$ rather
than to $\cf(\gk,{\gk^+}^m)$, even though $\card m=\gk$ implies that
${\gk^+}^m<\gk^+$.

If $o(\gk)<\gk^+$ then this consideration is not a problem, since in
that case $\cf$ and $\cf^m$ are equal at $\gk$.  This may be proved by
using iterated ultrapowers to compare the models $m$ and
$\ult(m,\cf(\gk,\gb))$ for $\gb<o^\cf(\gk)$, using the fact that for
$\gb<{\gk^+}^{(m)}$ the coherence functions $\cohere(\gk,\gb',\gb)$ of
$\coref$ are in $m$. In particular, this applies to the proof of
theorem~\mainthm(i) in section~4, where there is a fixed bound
$\gb_{0}$, smaller than the first relevant measurable  cardinal, for
the order $o(\ga)$ of any measurable cardinal $\ga $.  This bound also
makes unnecessary many of the complications of this section.  It means
that if $c\in \cc(\ga,\gb)$ then $\gb=o(c)$, and since it can be
assumed that $\gb_{0}$ is contained in $X$ there is no concern about
the definability of $\gb<\gb_{0}$. 
\medskip
Most of the rest of this section will be concerned with looking in
more detail at the structure of the covering set
$h^{m}\image(\gd;\cc^X)$. 
We will  generally  write  $m$  for  $m^{X}$,  and  in  general  will
frequently  drop  superscripts  when  they  are  not  necessary  to  prevent
ambiguity.

\defin A finite
increasing sequence $d=(d_{0},...,d_{n-1})$ of ordinals  is a {\it
weak support} in $\cc^{X}$ if for each $i<n$ there are
$\ga_{i},\gb_{i}\in h\image (d\restrict i)$ such that $d_{i}\in
\cc^{X}(\ga_{i},\gb_{i})$.
\endth

\lemma 
For every  set $y\in\range(\pi^*)$ there is a weak support $d$ such
that $y\in h\image(\rho^X\cup d)$.
\endth
\proof
Since every such set $y$ has the form $h(\nu)$ for some ordinal $\nu$,
it is sufficient to prove the lemma for the case where $y=\nu$ is an
ordinal.
The proof is by induction on $\nu$.  If $\nu$ is not an indiscernible
then $\nu=h(\nu')$ for some $\nu'<\nu$, and the induction hypothesis
implies that there is a sequence  $d\subset\nu'+1\subset\nu$ which is
a weak support for
$\nu'$ and hence for $\nu$.
If $\nu$ is an indiscernible then $\ga^m(\nu)$ and $\gb^m(\nu)$ are in
$h\image\nu$ and it follows by the induction hypothesis that they have
a weak support $d\subset\nu$.  Then $d\cup\{\nu\}$ is a support for $\nu$.
\endproof

For the lemma~\wksupportLem\ we need one more fact using the fine
structure:

\lemma\thmtag\hsupxiLem 
The skolem function $h$ for $m$ is the increasing union of a set of
partial functions $h^{\xi}$ such that each $h^\xi$ is a member of $m$.
Cofinally many of the functions $h^\xi$ are of the form $\pi^*(\bar
h^{\bar\xi})$ where $\bar h^{\bar\xi}$ is contained in $\bar h$, and
cofinally many of the functions $h^\xi$ are in $h^m\image\rho^X$.
\endth

\proof
The existence of the functions $h^\xi$ depends on the fine
structure of $m$.  We will give the  definition for the  case in which
$m$ is a $\gS_1$-code, that is, $m=J_\ga(\cf^m)$ for some ordinal 
$\ga$ and $h$ is the $\gS_1$ skolem
function for $J_\ga(\cf^m)$ with parameter $p=p_1^m$.  The
construction is similar for the general case of a $\gS^*_n$ code but 
does, of course, require knowledge of the
$\gS^*_n$-codes of a structure $J_\ga(\cf^m)$.

In the $\gS_1$ case the skolem function $h$ for $m$ is the partial
function defined in
$J_\ga(\cf^m)$ by
$$h(\nu)=x\iff\exists y\;
\bigl(R(\nu,x,p,y)\land\forall(x',y')<^m(x,y)\;\lnot R(\nu,x',p,y')\bigr)$$
where $\exists y\,R$ is the universal $\gS_1$ formula
and $<^m$ is the order of
construction of $J_\ga(\cf^m)$.
If $\ga$ is a limit ordinal and $\xi<\ga$ then we define $h^\xi$
for $\xi<\ga$ by $h^\xi(\nu)=x$ if and only if
$$x\in J_\xi(\cf^m)\land\exists y\in J_\xi(\cf^m)\;
\bigl(R(\nu,x,p,y)\land\forall(x',y')<^m(x,y)\;\lnot
R(\nu,x',p,y')\bigr).$$
If $\ga=\gg+1$ is
a successor ordinal then $h^\xi(\nu)$ is defined for integers $\xi<\gw$
by setting $h^\xi(\nu)=x$ if and only if
$$x\in S_{\gg\cdot\gw+\xi}(\cf^m)\land\exists y\in
S_{\gg\cdot\gw+\xi}(\barcf)\;
\bigl(R(\nu,x,p,y)\land\forall(x',y')<^m(x,y)\;\lnot
R(\nu,x',p,y')\bigr),$$
where $S_{\gg\cdot\gw+\xi}$ is the union of the images of
$J_{\gg}(\cf)\cup\{J_{\gg}(\cf^m)\}$
under the first $\xi$ of the $\cf^m$-rudimentary functions.
In either case it is clear that the functions $h^\xi$ form an increasing
sequence of partial functions, that each of the functions
$h^\xi$ is in $m$, and that the union of this sequence is the skolem
function $h$.

Let $\bar h^{\bar\xi}$ be defined in $\barm$ in the same way that
$h^\xi$ was defined in $m$. 
Since the skolem function $h$ for $m$ is the image under $\pi^*$
of the skolem function $\bar h$ for $\barm$, we have
$h^\xi=\pi^*(\bar h^{\bar\xi})$, as required in the first half of the
final sentence, whenever $\xi=\pi^*(\bar\xi)$.
If $\ga$ is a limit ordinal then $\range(\pi^*)$ is cofinal in $\ga$ since
$\pi^*$
is defined as an extender using functions which are members of
$\barm$. If $\ga$ is a successor ordinal then trivially every function
$h^\xi$ for $\xi<\gw$ is in the range of $\pi^*$. This proves the
first half of the last sentence of the lemma.

The second half of the final sentence uses a similar argument applied
to the function $i_{0,\gth}$ instead of to $\pi^*$.  If $\ga$ is a
limit ordinal then
the range of $i_{0,\gth}$ is cofinal in $\bar\ga$ since $i_{0,\gth}$
is an iterated ultrapower using functions in $\barm_0$, and hence
$\range(\pi^*\cdot i_{0,\gth})$ is cofinal in $\ga$.  Hence cofinally
many of the functions $h^\xi$ are in the range of
$\pi^*\cdot i_{0,\gth}$, and all such functions $h^\xi$ are in
$h\image\rho^X$. Again, if $\ga$ is a successor than all of the
functions $h^\xi$ for $\xi<\gw$ are in the range of
$\pi^*\cdot i_{0,\gth}$ and hence are in $h\image\rho^X$.
\endproof

\lemma\thmtag\wksupportLem  If $d$ is a weak support in
$\cc^{X}$,  $c\in \cc^{X}(\ga,\gb)$, and $d\cap [c,\ga)=\nothing$
then for all  $x\in h\image(c\cup d\cup \{\ga,\gb \})$ we have $c\in
x $ iff  $x\in \cf^{m}(\ga,\gb)$.
\endth
\proof
Suppose that the lemma is false, so that
there is a set $x\in h\image(c\cup d\cup\{\ga,\gb\})$
such that $c\in x\iff x\not\in\cf^m(\ga,\gb)$. We
begin by using the
last lemma to show that the  lemma is also false in
$\barm$ and we will then be able to reach a contradiction by using
the fact that the system
$\barcc$ of indiscernibles for $\barm$ is constructed from an
iterated ultrapower.

Let $h^\xi$ be a function as in lemma~\hsupxiLem\ such that
$h^\xi\in\range(\pi^*)$ and $x\in
h^\xi\image(c\cup d\cup\{\ga,\gb\})$.  Then $m$ satisfies the sentence
$$\exists x\in h^\xi\image(c\cup d\cup\{\ga,\gb\})\; \bigl(c\in x\iff
x\notin\cf^m(\ga,\gb)\bigr).$$ By elementarity, $\barm$ satisfies the
sentence $$\exists x\in \bar h^{\bar\xi}\image(\bar c\cup \bar
d\cup\{\bar\ga,\bar\gb\})\; (\bar c\in x\iff
x\notin\cf^{\barm}(\bar\ga,\bar\gb)).$$ This gives the required
counterexample in $\barm$. Dropping the bars, this means that we
have $c$, $\ga$, $\gb$, $d$ and $x$ in $\barm$ such that
$$
c\in\barcc(\ga,\gb),\qquad
d\cap[c,\ga)=\nothing,\qquad
x\in\bar h\image(c\cup d\cup\{\ga,\gb\})$$
and
$$
c\in x\iff x\notin\barcf(\ga,\gb).
$$

Now we claim that $x\in\bar h\image(c\cup d\cup\{\ga\})$, so that $\gb$
can be omitted.
If $\ga=\ga^{\barm}(c)$ then this is immediate since $\gb=\gb^{\barm}(c)\in
\bar h\image c$.
If $\ga\not=\ga^{\barm}(c)$ then $\ga^{\barm}(\ga)=\ga^{\barm}(c)$, and by the
construction of the system $\barcc$ of indiscernibles for $\barm$ we
have $\gb^{\barm}(\ga)\in \bar h\image c$ and
$\gb=\cohere^{\barm}\bigl(\ga^{\barm}(c),\gb',\gb^{\barm}(\ga)\bigr)(\ga)$ for
$\gb'=i_{\xi,\gth}(\gb)<\gb^{\barm}(\ga)$ in $\bar h\image c$. Since $\ga^{\barm}(c)$, $\gb'$,
and $\gb^{\barm}(c)$ are all in $\bar h\image c$ we have $\gb\in
\bar h\image(c\cup\{\ga\})$ and hence $x\in \bar h\image(c\cup d\cup\{\ga\})$.

Now  we claim that $d$ can also be omitted, so
that $x\in\bar h\image(c\cup\{\ga\})$. Every member $y$ of
an iterated ultrapower has a support which is minimal in the sense
that it  is contained in any support for $y$. Now if
$\ga=\gk_{\nu}=\crit(i_{\nu,\gth})$ then $x\subset\ga$ implies that
$x\in\range(i_{\nu+1,\gth})$
and hence $x$ has a support which is contained in $\ga+1$. On the
other hand $x\in \bar h\image(c\cup d\cup\{\ga\})$ implies that there is
a support contained in $c\cup d\cup\{\ga\}$, and
hence the minimal support must be contained in the intersection of
the sets $\ga+1$ and $c\cup d\cup\{\ga\}$.  
Since $d\cap[c,\ga)=\nothing$ this intersection is equal to 
$c\cup\{\ga\}$.

Now if $\ga=\ga^{\barm}(c)$ then it follows immediately that $c\in
x\iff x\in\cf^{\barm}(\ga,\gb)$. If $\ga<\ga^{\barm}(c)$ then
$x\in\cf^{\barm}(\ga,\gb)$ iff $i_{\nu,\gth}(x)
\in\cf^{\barm}(\ga^{\barm}(c),\gb^{\barm}(c))$ iff $c\in
i_{\nu,\gth}(x)$ iff $c\in x$ since $i_{\nu,\gth}(\ga)=\ga^{\barm}(c)$
and $i_{\nu,\gth}(\gb)=\gb^{\barm}(c)$.
\endproof

Our next goal is to show that $\cc^{X}$ is a system of indiscernibles for
$\coref$ rather than just for $m^{X}$ .  We begin with $\gw $ sequences, and
then use this special case to understand longer sequences. 

\defin \thmtag\indiscDef
\roster
\item"(i)" An  increasing  $\gw$-sequence  $\vc$  of  ordinals  is  an
{\it indiscernible  sequence}  (over  $\coref$) for  $(\ve\ga,\ve\gb)$  iff
for each function $g\in\coref$ there is $i_0\in\gw$ such that
$$\forall i>i_0\forall x\in
g\image(c_{i}\cup\set{\ga_{i},\gb_{i}})\quad(x\in\cf(\ga_{i},\gb_{i})\iff
c_{i}\in x)$$ 
\item"(ii)"
$\cd$ is a {\it system of indiscernibles} for $\coref$ if $\cd$ is a
function such that
\itemitem{(a)}
$\domain(\cd)\subset\domain(\cf)$, and
$\set{\cc(\ga,\gb):(\ga,\gb)\in\domain(\cd)}$ is a family of disjoint
subsets of $\ga$ for each ordinal $\ga$.
\itemitem{(b)}
 For all $(\ga,\gb)$ in $\domain(\cd)$, and for all $\nu \in \cd(\ga,\gb)$ and all $\gg <o(\nu)$ there is $\gg '<\gb $ and $\xi <\nu $ such
 that $\cd(\nu,\gg)\setminus\xi=\cd(\ga,\gg ')\cap [\xi,\nu)$. 
\itemitem{(c)}
If $c<\ga$ is a limit point of $\union_{\gb }\cd(\ga,\gb)$ of
cofinality $\gw $ then $c\in \cd(\ga,\gg)$ for some $\gg $ such that
$\union_{\gb \ge \gg }\cd(\ga,\gb)$ is bounded in $c$. 
\itemitem{(d)}
If $\vc$, $\ve\ga$ and $\ve\gb$ are sequences
such that $\vc$ is increasing, $c_{i}\in \cd(\ga_{i},\gb_{i})$ for
all $i\in \gw $ and $\union_{i}c_{i}=\union_{i}\ga_{i}$ then $\vc$ is an
indiscernible sequence for $(\ve\ga,\ve\gb)$. 
\endroster\endth

\lemma\thmtag\indiscLem (i) $\cc^{X}\restrict \gk $ is a system of indiscernibles for
$\coref$.  Furthermore clauses (a) and (b) of the definition of a system
of indiscernibles for $\coref$ hold for all of $\cc^{X}$, as do
clauses~(c) and~(d) if $o(\gk)<\gk^{+}$. 

(ii) Conversely, suppose that $\vc$ is an indiscernible sequence for
$(\ve\ga,\ve\gb)$, that $\vc$ and $\ve\ga $ are in $X$ and that either
$\ga_{i}<\gk$ for all 
$i$ or $o(\gk)<\gk^{+}$.  Then $c_{i}\in \cc^{X}(\ga_{i},\gb_{i})$
for all but finitely many $i\in \gw $.  \endth

\proof  We first prove that $\cc^{X}\restrict \gk $ is a system of
indiscernibles for $\coref$.  Clause~ii(a) and~ii(b) of
definition~\indiscDef\ are immediate for  $\cc^{X}$.  In order to
establish clause~ii(c), let $\vc$ be a strictly increasing sequence
with $c_{i}\in \cc^X(\ga,\gb_{i})$, and let  $c=\union_{i}c_{i}$. Note
that since $^{\gw }X\subset X$ we have $\vc\in X$ and hence $c\in X$.
We can assume that $\ga$ is not an indiscernible, since otherwise if
the assertion fails for $\ga$ then it still fails if  we replace $\ga$
with $\ga^m(\ga)$, using the fact that $\gb^m(\ga)\in h\image
c_i\subset h\image c$. It follows from the construction of $\barcc$
and the definition of $\cc^{X}$ that $c\in\cc^{X}(\ga,\gb)$ for some
$\gb<o(\ga)$. Suppose \iwlog\  that $\gb\le\gb_i$ for all $i<\gw$.

Let $\bar c$, $\bar c_i$, $\bar \gb$, and $\bar\gb_i$ be the preimages
under $\pi^*$ of $c$, $c_i$, $\gb$ and $\gb_i$ respectively. Then there are
$\nu_i<\gth$ such that $\bar c_i$ is the critical point of
$i_{\nu_i,\gth}:\barm_{\nu_i}\to\barm=\barm_\gth$, and $\bar c$ is the
critical point of $i_{\nu,\gth}$ where $\nu=\sup_i\nu_i$. Let
$\gd=o^{\barcf}(\bar c)$, so that $\barm_{\nu+1}=\ult(\barm_\nu,\barcf_\nu(\bar c,\gd))$ and $\bar\gb=i_{\nu,\gth}(\gd)$.
Then $\barcf_\nu(\bar
c,\gd)\notin N$, since otherwise it would be in $\barcore$ and hence
would not have been used in the iteration. We will reach a
contradiction by proving that $\barcf_\nu(\bar c_,\gd)$ is in $N$.

Suppose that $x\in\barm\cap\ps(\bar c)$.  Then
$x=i_{\nu_i,\nu}(x\cap\bar c_i)$ for all sufficiently large $i<\gw$,
and for these $i$ we have $x\in\cf_{\nu}(\bar c,\gd)$ iff $x\cap\bar
c_i\in\cf_{\nu_i}(\bar c_i,\gd_i)$ where $i_{\nu_i,\nu}(\gd_i)=\gd$.
If $\gb_i=\gb$ then it follows that $\bar c_i\in x$ if and only if $x\cap
\bar c_i\in\cf_{\nu_i}(c_i,\gd_i)$ if and only if
$x\in\cf^{\barm}(\bar c,\gd)$. Thus if $\gb_i=\gb$ for infinitely many
$i$ then $\cf_{\barm}(\bar c,\gd)=\set{x:\exists i_0\forall
i>i_0\,(\bar\gb_i=\bar\gb\implies c_i\in x}\in N$, so we can assume
\iwlog\
that $\gb<\gb_i$ for all $i$. Then $$i_{\nu_i,\gth}(\gd_i)=
i_{\nu,\gth}(\gd)=\bar\gb<\bar\gb_i=i_{\nu_i,\gth}(o^{\barm}(\bar
c_i)),$$  so $\barcf_{\nu_i}(\bar c_i,\gd_i)=\barcf(\bar c_i,\gd_i)\in
N$. Then $\left(\,\barcf(\bar c_i,\gd_i):i<\gw\,\right)\in N$ since
${^\gw}N\subset N$, and hence $\cf_{\nu}(\bar c,\gd)=\set{x\subset
\bar c:\exists
i_0\forall i>i_0\;x\cap \bar c_i\in\barcf(\bar c_i,\gd_i)}\in N$.  This
contradiction completes the proof of clause~ii(c) of
definition~\indiscDef.
\medskip

We now establish condition~ii(d) of definition~\indiscDef.  Suppose
that $c_{i}\in\cc^{X}(\ga_{i},\gb_{i})$ for almost all $i$ and
$\union_{i}c_{i}=\union_{i}\ga_{i}$ but $\vc$ is not an indiscernible
sequence for $(\ve\ga,\ve\gb)$.   The sequences $\vc$ and
$(\ve\ga,\ve\gb)$ are in $X$ since ${^\gw}X\subset X$, and by
elementarity the statement that $\vc$ is not an indiscernible sequence
for $(\ve\ga,\ve\gb)$ is true in $X$.  Let $g\in X$ be a witness, so
that for infinitely many i there is $x_{i}\in
g\image(c_{i}\cup\{\ga_{i},\gb_{i}\})$ such that $c_{i}\in x_{i}\iff
x_{i}\notin\cf(\ga_{i},\gb_{i})$.  We may assume \iwlog\ that this
equivalence holds for all $i$.  Since $g\in X$, there is a weak
support $d$ such that $g\in h\image d$, and hence $x_i\in
h\image(c_i\cup d\cup\{\ga_i,\gb_i\})$.  The assumption $\union
c_{i}=\union\ga_{i}$  implies that $d\cap [c_{i},\ga_{i})=\nothing$
for all but finitely many $i$, but then lemma~\wksupportLem\ implies
that $c_i\in x_i\iff x_i\in\cf(\ga_i,\gb_i)$ for all but finitely many
$i$. This contradiction completes the proof of clause~ii(d) of
definition~\indiscDef.

The argument above used the fact that
$\cf(\ga_{i},\gb_{i})=\cf^{m}(\ga_{i},\gb_{i})$ and hence is valid so
long as $\ga_{i}<\gk $ or $\ga_{i}=\gk $ and $o(\gk)<\gk^{+}$.
\medskip

We now prove clause~(ii) of lemma~\indiscLem.  Assume that $\vc$ is
an indiscernible  sequence for $(\ve\ga,\ve\gb)$ and that $\vc$ and
$\ve\ga $ are in $X$.  First we show  that $c_i$ is an indiscernible
and that $\ga^{m}(c_{i})=\ga^{m}(\ga_{i})$ for all but finitely many
$i$. (Recall that $\ga^{m}(\ga)=\ga $ if $\ga $ is not an
indiscernible.) Note that $\ga^{m}(c_{i})\ge\ga_{i}$ for all but
finitely many $i$ since
$c_i\le\ga^{m}(c_{i})\in h^{m}\image (c_{i})$ and $\vc$ is an indiscernible
sequence for $(\ve\ga,\ve\gb)$.  Since $c_{i}<\ga_{i}$ it follows
that $\ga^{m}(c_{i})\ge\ga^{m}(\ga_{i})$ for all but finitely many
$i\in\gw$. The other direction will
follow if we can show that $\ga^{m}(\ga_{i})\in h^{m}\image c_{i}$
for all but finitely many $\ga$.
To see this let $\nu_i$ be  the
least ordinal $\nu$ such that $\ga^m(\ga_i)\in h\image\nu$.
Then $\nu_{i}<\ga_i$ since $\ga^{m}(\ga_{i})\in h^m\image\ga_i$, and
since there is a function $g\in\coref$ such
that $\nu_i=g(\ga_i)$ for all $i$ the assumption
that $\vc$ is an indiscernible sequence for $(\ve\ga,\ve\gb)$ implies
that $\nu_i<c_i$ for all but finitely many $i\in\gw$.

Now we need to see that $c_{i}\in\cc^{X}(\ga_{i},\gb_{i})$ for almost
all $i$.  Let us first consider the special case in
which $\ve\gb\in X$.  We consider two subcases.  For the first
subcase we assume that $\ga^{m}(c_{i})=\ga_{i}$.  Then $\ga_i$ is not
an indiscernible and since $\gb_i\in X$ it follows that $\gb_{i}$ is
in $h^{m}\image\ga_{i}$.  It follows that $\gb_{i}$ is in
$h^{m}\image c_{i}$ for almost all $i$: otherwise define a function
$g$ by setting $g(\gg)$ equal to  the least ordinal $\nu $ such that
$\gg\in h^{m}\image\nu$.  Then the function $g$ is in $\coref$, but
$g(\gb_{i})$ is strictly  between $c_{i}$ and $\ga_{i}$ for
infinitely many $i$, contradicting the assumption that $\vc$ is an
indiscernible sequence for $(\ve\ga,\ve\gb)$.  Since
$\gb^{m}(c_{i})\in h^{m}\image c_{i}$ as well,
$\gb_{i}\not=\gb^{m}(c_{i})$ would imply that the least set $x$ in
$\cf(\ga_{i},\gb_{i})\setminus\cf(\ga_{i},\gb^{m}(c_{i}))$ is in
$h^{m}\image c_{i}$, so that with finitely many exeptions $c_{i}$
must  be a member of both $x$ and $\ga_{i}\setminus x$.  This
contradiction shows that
$\gb_{i}=\gb^{m}(c_{i})$, and $c_{i}\in\cc^{X}(\ga_{i},\gb_{i})$, as
required.  The other subcase,
$\ga^{m}(c_{i})=\ga^{m}(\ga_{i})>\ga_{i}$, is similar but slightly
more complicated.  There are $\gg_{i}$ such that
$\gb_{i}=\cohere^{m}(\ga^{m}(c_{i}),\gg_{i},\gb^{m}(\ga_{i}))(\ga_{i})$.
The ordinal $\gg_{i}$ is given from the ordinals $\gb_{i}$ and $\ga_{i}$ by a
function in $\coref$, and $\gg_{i}\in h^{m}\image\ga_{i}$, so the
least $\nu $ such that $\gg_{i}\in h^{m}\image\nu$ must be below
$c_{i}$ for almost all $i$.   If $\gb^{m}(c_{i})\not=\gg_{i}$
then the least $x$ such that
$x\in\cf^{m}(\ga^{m}(c_{i}),\gb^{m}(c_{i}))\setminus\cf^{m}(\ga^{m}(c_{i}),\gg_{i})$
is in $h^{m}\image c_{i}$ and so $c_{i}$ is again in both $x$ and
$\ga_{i}\setminus x$ for almost all $i$. The contradiction shows that
$\gb^{m}(c_{i})=\gg_{i}$.
To complete the proof that
$c_i\in\cc^X(\ga_i,\gb_i)$ we must verify that $\gb^m(\ga_i)\in h^m\image
c_i$ for almost all $i$. Now $\gb^m(\ga_i)\in h^m\image(\ga_i)$, so if $g(\xi)$ is
the least ordinal  $\nu$ such that $\xi\in h^m\image\nu$ then $g$ is
in $\coref$  and $g(\gb^m(\ga_i))<\ga_i$ for all $i$.
Since $\gb^m$ is also in $\coref$ it follows that
$g(\gb^m(\ga_i))<c_i$ for
almost all $i$. This completes the proof for the case when $\vec\gb$
is in $X$.

Now let $\ve\gb$ be an arbitrary sequence such that $\vc$ is an
indiscernible sequence for $(\ve\ga,\ve\gb)$.  By elementarity there
must be a sequence $\ve\gb'$ in $X$ such that $\vc$ is an
indiscernible sequence for $(\vec\ga,\vec\gb')$, and by the last
paragraph $c_{i}\in\cc^{X}(\ga_{i},\gb'_i)$ for almost all $i$.  If
$\gb_{i}\not=\gb_i'$ for infinitely many $i$ then $\vc$ is an
indiscernible sequence for two different sequences $\ve\gb$ and
$\ve\gb'$, and by elementarity there is a second sequence $\gb''$ in
$X$ so that $\vc$ is an indiscernible sequence for
$(\ve\ga,\ve\gb'')$ and $\gb''_i\not=\gb'_i$ for infinitely many $i$.
This is impossible since, again by the last paragraph, we would also
have $c_{i}\in\cc^{X}(\ga_{i},\gb'')$ for almost all $i$.\endproof

\cor\thmtag\cXeqcXpr Suppose that $\vc$, $\ve\ga$ and $\ve\gb$ are sequences such that
$c_{i}\in\cc^{X}(\ga_{i},\gb_{i})$ for all $i\in\gw $,
$\union_{i}c_{i}=\union_{i}\ga_{i}$, and either $\ga_{i}<\gk $ or
$o(\gk)<\gk^{+}$.  Then for any $X'\supset X$ we have $c_{i}\in\cc^{X'}(\ga_{i},\gb_{i})$ for all but finitely many $i$.\endth

Proof: By applying lemma~\indiscLem(i) to $\cc^{X}$ we get that $\vc$
is an indiscernible sequence for $(\vec\ga,\vec\gb)$, and by applying
lemma~\indiscLem(ii) to $\cc^{X'}$ we then get that
$c_{i}\in\cc^{X'}(\ga_{i},\gb_{i})$ for all but finitely many $i$.
%%%\endproof

The next corollary may be regarded as a semicontinuity property for the
function $\gb(c)$.  Note that it is trivial if $X$ is not cofinal in $c$.  We
conjecture that it is false if ``bounded in $c$'' is replaced by ``bounded in
$c\cap X$''.

\cor \thmtag\semicont
For all but finitely many indiscernibles $c$, if $\ga$ and $\gb$ are
any ordinals such that $c\in\cc(\ga,\gb)$ then
$\union_{\gb '\ge\gb}\cc^{X}(\ga,\gb')$ is bounded in $c$.
Thus, except at these
finitely many exceptional points, $$\gb(c)\ge\limsup_{c'\to
c}\set{\gb(c')+1:c'<c\text{ and }\ga (c')=\ga(c)}$$ where the
$\limsup$ is taken to be $0$ if the set of relevant ordinals $c'$ is
bounded in $c$.
\endth

\proof 
Suppose that the lemma is false, and let
 $\vc$, $\ve\ga $, and $\ve\gb$ be $\gw$-sequences such that
that for each $i\in\gw$ we have
$c_{i}\in\cc^{X}(\ga_{i},\gb_{i})$ but $\union_{\eta\ge\gb_i
}\cc^{X}(\ga_{i},\eta)$ is unbounded in $c_{i}$.  By
lemma~\indiscLem(i) and definition~\indiscDef(iic) each $c_{i}$ must have
uncountable cofinality, and we may assume that $\vec c$ is strictly increasing.

We first show that we can assume that $\union_i c_i=\union_i \ga_i$.
This will involve two stages, the first of which is to show that we
can assume that $\ga_i=\ga(c_i)$ and $\gb_i=\gb(c_i)$.
Since $c_i\in\cc(\ga_i,\gb_i)$ there is $\gd_i<c_i$ such that
$\gb(\ga_i)$ is in
$h\image\gd_i$. Then
$\union_{\gl\ge\gb_i}\cc(\ga_i,\gl)\subset
\union_{\gl\ge\gb(c_i)}\cc(\ga(c_i),\gl)$,
since if $\gl'\in X\cap o(\ga_i)$ then there is $\gl<\gb(\ga_i)$ such
that $\gl'=\cohere(\ga,\gl,\gb(\ga_i))(\ga_i)$ and
$\cc(\ga_i,\gl')=\cc(\ga(c_i),\gl)\setminus\xi$ for some $\xi<\ga_i$.
It follows that the corollary is also false for $\ga_i=\ga(c_i)$ and
$\gb_i=\gb(c_i)$.

If there is an infinite subset
$I\subset\gw $ such that $\ga (c_{i})<\ga (c_{i'})$ for $i<i'$ in $I$ then
$c_{i}\le\ga (c_{i})<c_{i'}$ for $i<i'$ in $I$ since
$\ga(c_i)\in h\image c_i$, and   hence $\union_{i\in
I}\ga(c_{i})=\union_{i\in I}c_{i}$. Thus we can assume that there is no
such set $I$
and hence that there is an ordinal $\ga$ such that $\ga(c_i)=\ga$ for
all sufficiently large $i<\gw$.
Set $c=\union_{i\in\gw }c_{i}$.  Then $\ga(c)\le\ga$ since
$\ga=\ga(c_i)\in h\image c_i\subset h\image c$ for all $i$, and
$\ga\le\ga(c)$ since $\ga(c)\in h\image c_i$
for all sufficiently large $i\in\gw$.  Thus we can assume \iwlog\ that
$\ga (c)=\ga (c_{i})=\ga$ for all $i$.
Since $\cof(c)=\gw$ the corollary is true for $c$ and hence
$\gb(c)>\gb(c_i)$
for all sufficiently large $i<\gw$.
Assume \iwlog\ that this is true for all $i$, and in addition there is
$\xi<c_0$ such that
$\gb(c)\in h\image\xi$ and
$\union_{\gb'\ge\gb(c)}\cc^X(\ga,\gb')\cap c$ is contained in $\xi$.
Now define new sequences $\vec\ga'$ and $\vec\gb'$ by setting
$\ga'_i=c$ and
$\gb_i'=\cohere(\ga,\gb_i,\gb(c))(c)$. If $\xi<\nu<c_i$ and
$\nu\in\cc(\ga,\eta)$ where $\eta\ge\gb_i$ then $\eta<\gb^m(c)$ and
$\nu\in\cc(c,\eta')$
where $\eta'=\cohere(\ga(c),\eta,\gb(c))(c)>\gb'_i$.
Hence $\union_{\eta\ge\gb'_i}\cc^{X}(\ga'_{i},\eta)\cap c_i=
\union_{\eta\ge\gb_i}\cc^{X}(\ga_{i},\eta)\cap c_i$ is 
unbounded in $c_i$ for all $i\in\gw$. Thus the sequences $\vc$, $\vec\ga'$ and
$\vec\gb'$ form a witness to the failure of the corollary such that
$\union_i c_i=\union_i\ga'_i$.

Thus we can assume that $\union_i\ga_i=\union_i c_i$.  For each
$i\in\gw$ we define $U_{i}$ to be the set of all subsets
$x$ of $c_{i}$ such that for all sufficiently large $d\in
c_i\cap\union_{\gb'\ge\gb_i}\cc(\ga_{i},\gb')$ we have 
$$\alignat2
d&\in x&&\qquad\text{if $d\in\cc(\ga_i,\gb_i)$}\\ 
x\cap
d&\in\cc\bigl(c_i,\cohere(\ga_i,\gb_i,\gb')(c_i)\bigr)&&\qquad\text{if
$d\in\cc(\ga_i,\gb')$ with $\gb_i<\gb'$.}
\endalignat
$$ The filter $U_{i}$ is countably complete because $\cof(c_{i})>\gw$.
It is also normal: otherwise let $f\in\coref$ be such that
$\set{\nu:f(\nu)<\nu}\in U_i$. If $d\in\cc(\ga_i,\gb_i)$ then set
$f'(d)=f(d)$, and if $d\in\cc(\ga_i,\gb')$ where $\gb'>\gb_i$ then set
$f'(d)=\gd$ where
$\set{\nu<d:f(\nu)=\gd}\in\cf(d,\cohere(\ga_i,\gb_i,\gb')(d))$. If
$f'$ is eventually constant, say $f'(d)=\gd$ for sufficiently large
$\nu$, then $\set{\nu:f(\nu)=\gd}\in U_i$, was was to be shown. If
$f'$ is not eventually constant then there is an infinite increasing
sequence $\vec d$ such that $f'(d_i)<d_i$ for all $i$ and $f'(d_i)$ is
not eventually constant. Let $d=\union_id_i$ so that
$\gb^m(d)\ge\limsup_i\gb(d_i)$ by definition~\indiscDef(iic). If $U$ is
defined on $d$ as $U_i$ was defined on $c_i$ then
$U=\cf(d,\cohere(\ga_i,\gb_i,\gg)(d))$ where $d\in\cc(\ga_i,\gg)$, but
this is a contradiction because $f$ witnesses that $U$ is not
normal.

Thus $U_{i}=\cf(c_{i},\gd_{i})$ for some $\gd_{i}<o(c_{i})$.  Now take
$X'\supset X$ such that $\cc^{X}\in X'$.  Then $\gd_{i}\in X'$ for
all $i$ because $U_{i}\in X'$.  For all but finitely many $i$ we have
$c_{i}\in\cc^{X'}(\ga_{i},\gb_{i})$  and hence
$\gd_{i}=\cohere(\ga_{i},\gb'_{i},\gb_{i})(c_{i})$ for some
$\gb'_{i}\in h^{X'}\image (c_{i}\cup \{\ga_{i}\})\cap\gb_{i}$. Now
let $x_{i}$ be the least set in $\coref$ which is in
$\cf(\ga_{i},\gb_{i})\setminus\cf(\ga_{i},\gb'_{i})$.  Then $x_{i}\in
h^{X'}\image(c_{i}\cup\{\ga_{i}\})$.  It follows that $x_{i}\cap
c_{i}\in U_{i}$ by the definition of $U_{i}$, but $x_{i}\cap
c_{i}\notin\cf(c_{i},\gd_{i})$ because
$x_{i}\notin\cf(\ga,\gb'_{i})$.  This contradiction completes the
proof of the corollary.\endproof

Now we know that $\cc^{X}$ is
a system of indiscernibles for $\coref$ rather than merely for
$m^{X}$ and that $w\subset X\subset
h\image(\rho^X\cup\union_{\ga,\gb}\cc(\ga,\gb))$.
In order to prove the SCH we would like to show that we need not
consider every system of the form $\cc^X$, but rather that there is a
small standard set of systems such that every set $w$ can be covered
by a system in the standard set. 
To do this we will need to come to a more precise understanding of the
similarity between different systems $\cc^X$ and $\cc^{X'}$ of
indiscernibles
than was given by corollary~\cXeqcXpr.
In the process we will
replace
$h\image(\rho^X\cup\union_{\ga,\gb}\cc(\ga,\gb))$ with a more delicate
covering: we will define functions $s^X$ and $a^X$ and then define
$h\image(\gd;\cc^X)$ to be the smallest set containing $\gd$ and closed
under the functions $s^X$ and $a^X$. Then we will show that
$h\image(\gd\cup\union_{\ga,\gb}\cc^X(\ga,\gb))=h\image(\gd;\cc^X)$ for
all $\gd$, so that $w\subset h\image(\gd;\cc^X)$. Finally we will show
that if $X$ and $X'$ are two different sets then $s^X$ and $a^X$
only differ from $s^{X'}$ and $a^{X'}$ on a bounded subset of $X\cap
X'$. In the next two sections this fact will be used to select the set
of standard systems $\cc^X$.  In section~4 we will actually use a
modification of the present definition of $h(\gd;\cc)$, and we will be
considering systems $\cc$ such that it is not true that there is
ordinal $\gd$ such that every indiscernible in $\cc$ is in
$h\image(\gd;\cc)$. Notice that our definition ensures that
$\card{h\image(\gd;\cc)}=\card\gd$, regardless of the size of $\cc$.

The first of the two functions used to define $h\image(\gd,\cc)$
is the least  indiscernible
function: $s^{\cc}(\ga,\gb,\gg)$ is the least member of $\cc(\ga,\gb)$
above $\gg$ (this definition is slightly modified below). We will
observe
that  under the assumption that $\set{o(\ga):\ga<\gk}$ is bounded
below $\gk$---that is, under the assumption of theorem~\mainthm\ for
countable cofinality---it is sufficient to close $h\image(\rho;\cc)$
under the functions $h$ and $s^{\cc}$. For larger sequences we need to
introduce the notion of accumulation points; the second function, $a^X$, used
in the definition of $h\image(\gd;\cc^X)$
is the least accumulation point function.

It should be noted that the definition of an accumulation point
depends directly on the set $X$ as well as on the system $\cc^X$
of indiscernibles. In particular, the domain of $\cc^X$ is contained
in that of $\cf^m$, while the ordinals $\gb$ in the definition of an
accumulation point are in
$X$ and hence come from the domain of $\cf$.  Thus the definition of
accumulation point for $(\ga,\gb)$ only makes sense if either $\ga<\gk$ or
$o(\gk)<\gk^+$, so that $\cf^m$ and $\cf$ agree at $\ga$.

\defin \thmtag\accptDef
(1) An ordinal $c$ is an {\it accumulation point} 
for $(\ga,\gb)$ in $X$ if  $c$,
$\ga$ and $\gb$ are in $X$, $\gb\le o(\ga)$, and for all $\gb '\in\gb\cap X$ the set
$$\union\set{\cc^{X}(\ga,\gb''):\gb''\in X \text{ and }\gb'\le\gb''<\gb }$$ is 
cofinal in $X\cap c$.  If $\ga=\ga(c)$ then this may be stated
equivalently
$$\limsup_{d\to c'}
\set{\gb(d):\ga (d)=\ga\text{ and }\gb (d)<\gb }=\sup(\gb\cap X).
$$
where $c'=\sup(c\cap X)$.

(2) We write $s$ for the least indiscernible function: 
$c=s(\ga,\gb,\gg)$ iff $\gg<c\in\cc(\ga,\gb)$ and $\union_{\gb
'\ge\gb }\cc(\ga,\gb ')\cap c\subset\gg +1$.

(3) We write $a$ for the least accumulation point function:
 $c=a(\ga,\gb,\gg)$ if $c>\gg $, $c$ is the least accumulation
point for $(\ga,\gb)$ above $\gg$, and $\union_{\gb'\ge\gb }\cc(\ga,\gb')\cap
c\subset\gg +1$.   

(4) A finite  sequence $d$ of ordinals is a (full) support in $X$ if
for each $i<\operatorname{length}(d)$  there is $(\ga,\gb)\in
h\image(d\restrict i)$ such that either $d_{i}=s(\ga,\gb,d_{i-1})$
or $d_{i}=a(\ga,\gb,d_{i-1})$, where for $i=0$ we take $d_{-1}=0$.
The sequence $d$ is a support for a set $e$ of ordinals if $d$ is a
support and every member $\nu$ of $e$ is either in $d$ or in
$h\image(d\cap\nu)$.
\endth

\lemma\thmtag\supportOne Assume that $c\in\cc(\ga,\gb)$ and either $\ga<\gk$ or $o(\gk)<\gk^+$.
Then there is  $\gg\in c\cap X$  such that either $c=s(\ga,\gb,\gg)$ or 
there is an $\eta$ in $h\image(\{\ga\}\cup c\cap X)$   
such that $c=a(\ga,\eta,\gg)$.
\endth
\proof  Assume that  $c\in\cc(\ga,\gb)$ and $c\not=s(\ga,\gb,\gg)$ 
for any $\gg\in X\cap c$,
and let $\eta\le o(\ga)$ be the least member of $X$ such that
$\union_{\gl\ge\eta}\cc(\ga,\gl)$ is bounded in $X\cap c$.
Then  $c$ is an accumulation point in $X$ for $(\ga,\eta)$.  First we will
show 
that there is $\gg\in X\cap c$ such that $c=a(\ga,\eta,\gg)$, and then we will
show that $\eta\in h\image (\{\ga\}\cup c\cap X)$.

Define an $\gw$~sequence $\ve \ga=\seq{\ga_{i}:i\in\gw}$ of
ordinals as follows: Set
$\nu=\union_{\gl\ge\eta}\cc^{X}(\ga,\gl)\cap c$, so that $\nu<c$
and $X\cap(c\setminus\nu)\not=\nothing$ by the choice of $\eta$.
Now pick $\ga_{0}\in X\cap (c\setminus\nu)$
and for $i>0$ define $\ga_{i}=a(\ga,\eta,\ga_{i-1})$.
If $\ga_i=c$ for some $i$ then we can take $\gg=\ga_{i-1}$, 
since $\ga_{i'}\in X$ for each $i'<i$.
Otherwise set $\ga'=\union_{i}\ga_{i}\le c$.  Then definition~\indiscDef ii(c)
implies that
$\ga'\in\cc^{X}(\ga,\tau)$ for some $\tau\ge\eta$.  If $\ga'<c$ then
this contradicts the choice
of $\ga_{0}$, while if $\ga'=c$ then $\tau=\gb$ and $c=s(\ga,\gb,\ga_{0})$,
contrary to assumption.

Now we have to show that $\eta\in h\image(\{\ga\}\cup X\cap c)$. As a
first step we will show that there is an ordinal $\nu\in X\cap c$
such that $h\image(\{\ga\}\cup \nu)$ is cofinal in $\eta\cap X$.   We
define an infinite sequence of indiscernibles: Pick $\nu_0\in X\cap
c$ so that $\nu_0\ge\sup(c\cap\union_{\gl\ge\eta}\cc(\ga,\gl))$. Now
suppose that $\nu_i$ is defined.  If $h\image\{\ga\}\cup\nu_{i}$ is
cofinal in $X\cap\eta$ then $\nu_i$ is the desired ordinal $\nu$;
otherwise let $\nu_{i+1}$ be the least ordinal such that
$\nu_{i+1}\in\cc(\ga,\gl_{i+1})$ for some $\gl_{i+1}\ge
\sup(h\image\{\ga\}\cup\nu_i)$. Then $\nu_{i+1}<c$ by the choice
of $\eta$. This process must stop in finitely many steps, for otherwise let
$\nu=\union_i\nu_i$.  Then corollary~\semicont\ implies that
$\nu\in\cc(\ga,\gl)$ for some $\gl$ such that $\gl>\gl_i$ for all
sufficiently large $i<\gw$ and since $\gl\in h\image(\{\ga\}\cup\nu_i)$ for
some $i<\gw$ this  contradicts the definition of $\nu_{i+1}$.

Now lemma~\hsupxiLem\ implies that there is a partial function
$h^\xi\subset h$ such that $h^\xi\in h\image\rho^X$ and $\eta\in
h^\xi\image(\ga+1)$. Let $\eta'$ be the least member of
$h\image(\{\ga\}\cup X\cap\nu)$ above $\eta$ and set
$f(\gg)=h^\xi(\ga,\gg)$ if $h^\xi(\ga,\gg)<\eta'$ and $f(\gg)=0$
otherwise.  Then $f\in h\image(\{\ga\}\cup(X\cap\nu'))$ for some
$\nu'<c$. We claim that $f\image c$ is cofinal in $\eta\cap X$: If not
then there is an ordinal $\gz\in h\image(\{\ga\}\cup c\cap X)$ such
that $\gz<\eta$ and $f\image c\subset\gz$.  Then the least ordinal
$\gd$ such that $f(\gd)>\gz$ is in $h\image(\{\ga\}\cup c\cap X)$, but
this is impossible because $c\le\gd<\ga$ and $c$ is an indiscernible
for a measure on $\ga$.

Finally, we claim that there is $\nu'\in c\cap X$ such that
$f\image\nu'$ is cofinal in $X\cap\eta$.  This will complete the proof
of the lemma, since then $\eta=\sup f\image\nu'$ is in
$h\image(\{\ga\}\cup X\cap c)$.
Define $g:\nu\to c$ by setting $g(\gz)$ equal to the least ordinal
$\gg<c$ such that 
$f(\gg)\ge h(\ga,\gz)$, and set $\nu'=\sup\range(g)$.  
Then $\nu'<c$ since $g$ witnesses that $\cof(\nu')=\nu$ while $c$ is regular
in $\coref$, and $f\image\nu'$ is cofinal in $X\cap\eta$ since
$h\image(\nu\cup\{\ga\})$ is cofinal in $X\cap\eta$.
\endproof

The next lemma will be used directly in section~4.  For section~3 we
will have to extend the proof to deal with the case $o(\gk)\ge\gk^+$,
where accumulation points are not defined at $\gk$.

\lemma\thmtag\havesupports If $o(\gk)<\gk^+$ then every finite set $e$ of
ordinals in $X\cap \coref\cap H_{\gk }$ has a support in $X$.\endth
\proof 
We will prove the lemma by 
induction on $\gz=\sup e$. If $\gz$ is not an indiscernible then there is
$\gg\in X\cap\gz$ such that $\gz=h(\gg)$.  By the induction hypothesis there
is a support $d$ for $\{\gg\}\cup e\cap\gz$, and then $d$ is also a support
for $e$. If $\gz$ is an indiscernible then
by lemma~\backthm1\
there are $\ga $ and $\gb $ in $h\image\gz$ and and $\gg\in X\cap\gz$
so that either $\gz=s(\ga,\gb,\gg)$ or $\gz=a(\ga,\gb,\gg)$.  By
the induction hypothesis there is a support $d$ for $e\cap\gz\cup\{\gg\}$ such
that
$\ga$ and $\gb$ are in $h\image d$.  Then $d$ is a support for $e$.
\endproof

Gitik \cite{G?} has shown that this lemma cannot be strengthened by
removing the accumulation point function, $a(\ga,\gb,\gg)$, from the
definition of a support. Gitik's construction requires a model in
which there is a measurable cardinal $\gk$ such that
$\set{o(\ga):\ga<\gk}$ is unbounded in $\gk$. It is known that a
cardinal such that $\set{o(\ga):\ga<\gk}$ is unbounded in $\gk$ is
required for the existence of accumulation points, but it is not known
whether the measurability of $\gk$ can be eliminated from the
hypothesis to Gitik's result, nor is it known whether a limit of
accumulation points can be singular in $\coref$. An affirmative answer
to the following problem would answer both of these questions
affirimatively.

\proc{Problem}Suppose that $\ga=(\ga_{i}:i\in\gw)$ is an increasing sequence
of measurable cardinals such that $o(\ga_{i+1})=\ga_{i}$.  Is there a larger
model $M$ in which each $\ga_{i}$ is still measurable and such that if $\ve\gg $
and $\ve\gb $ are any sequences such that $\gg_{i}<\ga_{i}$ and
$\gb_{i}<o(\ga_{i})$ for all $i\in\gw $, then there is a sequence $\vc$ which is
an indiscernible sequence for $(\ve\ga,\ve\gb)$ such that $\gg_{i}<c_{i}$
for all $i\in\gw$?\endth

If the answer is yes, then work in $L(\Cal U)$, with the ordinals $\ga_{i}$ as
in the hypothesis, and let $M\supset L(\Cal U)$ be the model asked for.  Let
$U_{i}$ be a measure on $\ga_{i}$ in $M$, and let $j:M@>>> N$ be the iterated
ultrapower obtained by using each of the measures $U_{i}$ once.  The model
$N[\ve\ga]$ is a generic extension of $N$ by a variant of Prikry forcing, and
hence has the same core model as $N$.  Let $\gb_{i}$ be the unique ordinal
$\gb <o^{\Cal U }(\ga_{i})$ such that $U_{i}\supset\Cal U (\ga_{i},\gb)$,
and in $N[\ve\ga]$ take $X\prec H_{\gk^{++}}$, where $\gk=\union_{i}\ga_{i}$,
with $\ve\ga$, $j(\ve\ga)$, and $j(\ve\gb)$ in $X$.  Since $\ve\ga$ is an
indiscernible sequence for $(j(\ve\ga),j(\ve\gb))$ we have
$\ga_{i}\in\cc^{X}(j(\ga_{i}),j(\gb_{i}))$ for all sufficiently large $i<\gw$. 
We claim that there is no sequence $\ve\gg $ such that
$\ga_i=s(j(\ga_{i}),j(\gb_{i}),\gg_{i})$ for all sufficiently large $i$. 
Suppose that $\ve\gg $ is such a sequence.  Then the model $M$ contains a
indiscernible sequence $\vc$ for $(\ve\ga,\ve\gb)$ such that $\gg_{i}<c_{i}$
for all $i\in\gw $.  It follows that $j(\vc)$ is an indiscernible sequence for
$(j(\ve\ga),j(\ve\gb))$ such that for all $i$ we have $\gg_{i}\le
j(\gg_{i})<j(c_{i})<\ga_{i}$ and by elementarity $X$ also satisfies that there
is such a sequence $\vd$.  Thus $s(j(\ga_{i}),j(\gb_{i}),\gg_{i})\le
d_i<\ga_i$ for all sufficiently large $i<\gw$. 

\medskip
We  now  end  section~2  by  showing  that  the  least  indiscernible
function  and  least  accumulation  point  function  are,  to  some  degree,
uniquely  determined.

\lemma  \thmtag\unique Suppose that $\ve\gg$, $\ve\ga$ and $\ve\gb$ are $\gw$ sequences
of ordinals contained in $X\cap X'$, that $\ve\gg$ is strictly
increasing, that $\union_{i}\gg_{i}=\union_{i}\ga_{i}$, and that
either $\ga_{i}<\gk $ or $o(\gk)<\gk^{+}$.  Then the following
equations hold for all but finitely many $i\in\gw$ 
$$
\alignat1
s^{X}(\ga_{i},\gb_{i},\gg_{i})&=s^{X'}(\ga_{i},\gb_{i},\gg_{i})\tag1\\
a^{X}(\ga_{i},\gb_{i},\gg_{i})&=a^{X'}(\ga_{i},\gb_{i},\gg_{i})\tag2
\endalignat$$
where the equals sign means that  if either side exists then both sides
exist and are equal.
\endth

\proof Consider equation~(1).  We can assume \iwlog\ that 
$c_{i}=s^{X}(\ga_{i},\gb_{i},\gg_{i})$ exists
for all $i\in\gw$.  Let $\phi(\ve c,\ve\ga,\ve\gb,\ve\gg)$ be a formula
asserting that 

\smallskip
{\narrower\narrower\noindent
The sequence $\ve c$ is an indiscernible sequence for $(\ve\ga,\ve\gb)$ such
that $c_i>\gg_i$ for all $i<\gw$, and there is no infinite subset
$I\subset\gw$ and sequences $\ve e=\seq{e_i:i\in I}$ and
$\ve\gb'=\seq{\gb'_i:i\in I}$ such that $\gg_i<e_i<c_i$ and $\gb_i\le\gb'_i$
for each $i\in I$ and $\ve e$ is an indiscernible sequence for
$(\ve\ga,\ve\gb')$.\par
}
\smallskip
Then  $\phi(\vc,\ve\ga,\ve\gb,\ve\gg)$ is true in $X$ and hence in $V$. Thus
$X'$ also
satisfies that there is a sequence $\vd$ satisfying 
$\phi(\vd,\ve\ga,\ve\gb,\ve\gg)$.  But then 
$\phi(\vd,\ve\ga,\ve\gb,\ve\gg)$ is also true in $V$, and 
this implies that $\vc$ and $\vd$ are eventually equal.

Now consider equation~(2).  Assume that $c_{i}=a^{X}(\ga_{i},\gb_{i},\gg_{i})$
exists for all $i\in\gw$.
Let $\phi(\vd,\ve\ga,\ve\gb)$ be a formula asserting that $\vd$ is a sequence
of accumulation points for $(\ve\ga,\ve\gb)$, that is (using
lemma~\indiscLem),

\smallskip
{\narrower\narrower\noindent
For any sequences $\ve\nu$ and $\ve\gl$ such that $\nu_{i}<d_{i}$ and
$\gl_{i}<\gb_i$ for all $i\in\gw$ 
there are sequences $\ve e$ and $\ve\gb'$ such that $\nu_i<e_i<d_{i}$,
$\gl_{i}\le\gb'_{i'}<\gb_{i}$, and $\ve e$ is an indiscernible sequence for
$(\ve\ga,\ve\gb')$. 
\par}
\smallskip\noindent
Then $X$ satisfies the formula $\theta(\vc,\ve\ga,\ve\gb,\ve\gg)$:
\smallskip
{\narrower\narrower\noindent $\phi(\vc,\ve\ga,\ve\gb)$ is true, but
$\phi(\ve e,\ve\ga,\ve\gb')$ is false for  all pairs $\ve e$ and
$\ve\gb'$ of sequences defined on an infinite subset $I$ of $\gw$ such
that  $\gg_i<e_i<c_i$ and $\gb_i'>\gb_i$ for all $i\in I$, and there
are no sequences $\ve e$ and $\ve\gb'$ defined on an infinite subset
$I\subset\gw$ such that  $\ve e$ is an indiscernible sequence for
$(\ve\ga,\ve\gb')$ and for all $i\in I$ we have $\gg_i<e_i\le c_i$ and
$\gb_i\le\gb'_i$.\par}
\smallskip

Then as in the argument for equation~(1), $V$ also satisfies
$\theta(\vc,\ve\ga,\ve\gb,\ve\gg)$ and hence $X'$ satisfies that there is a
sequence $\vd$ such that $\theta(\vd,\ve\ga,\ve\gb,\ve\gg)$.  Then
$\theta(\vd,\ve\ga,\ve\gb,\ve\gg)$ is also true in $V$, and it follows that
$\vc$ and $\vd$ must be equal except on an initial segment.
\endproof 

\cor \thmtag\uniqueFinally 
For any $X$ and $X'$ there is an ordinal $\xi<\gk$ such that
whenever $\ga,\gb,\gg\in X\cap X'$, $\gg >\xi$, and $\ga<\gk$
we have $s^{X}(\ga,\gb,\gg)=s^{X'}(\ga,\gb,\gg)$ and
$a^{X}(\ga,\gb,\gg)=a^{X'}(\ga,\gb,\gg)$.\endth

\proof If the corollary is false then there would be sequences
$\ve\ga$, 
$\ve\gb$ and $\ve\gg$ such that for each $i\in\gw$ we have 
$\gg_i<\ga_i\le\gg_{i+1}$ and the corollary is false for $\ga_i$, $\gb_i$ and
$\gg_i$, contradicting 
lemma~\unique.
\endproof

\heading{\newsectno3 Uncountable Cofinality}\endheading

  In this section we prove theorem~\mainthm(ii), which is the easier part of
the main theorem.  Most of the necessary tools are in section~2; the only
difficulty is the indiscernibles for measures on $\gk $, which are not covered
by lemma~\supportOne\ and corollary~\uniqueFinally.  We assume throughout this
section that $\gk $ is a singular strong limit cardinal of uncountable
cofinality. 

For each function $h:\gk@>>>H_{\gk}$ in $\coref$ and each ordinal
$\ga<\gk$ pick a set $X=X^{h,\ga }\prec H_{\gk^{++}}$ with $\card
X<\gk$ such that $\ga\cup\{h,\ga\}\subset X$ and ${^\gw}X\subset X$.
The main lemma of this section is

\lemma\thmtag\mainthree  For every set $w\subset\gk $ of cardinality less
than $\gk$ there is a function $h$ and a sequence $\ve\ga=(\ga_{k}:k\in\gw)$
such that $w\subset\union_{k\in\gw}X^{h,\ga_{k}}$. 
\endth

  Theorem~\mainthm(ii) follows easily from Lemma~\mainthree: Since $\coref$
satisfies the GCH there are only $\gk^{+}$ many functions $h$, and since $\gk$
is a strong limit cardinal of uncountable cofinality there are only $\gk$ many
sequences $\ve\ga$.  Hence there are only $\gk^{+} $ many pairs
$(h,\ve\ga)$ and hence only $\gk^+$ many sets $\union_{k\in\gw}X^{h,\ga_{k}}$. 
Since each of these sets has cardinality less than $\gk$ and $\gk$ is a strong
limit cardinal, each of these sets has
fewer than $\gk $ subsets and thus there are only $\gk^+$ subsets of $\gk$ of
cardinality less than $\gk$, and hence $2^\gk=\gk^+$.

We need some preliminaries before proving lemma~\mainthree.  Let
$X\prec H_{\gk^{++}}$ be arbitrary.  For $\nu<\gk$ define
$A^{X}_{\nu} $ to be  $\set{s(\gk,\gb,\nu):\gb <o^{m^X}(\gk)}$.  Then
$A^X_\nu$ is $\gw$-closed.  To see this, suppose that $c$ is a limit
point of $A^X_\nu$  of cofinality $\gw$, say $c=\union_i c_i$ where
$c_i=s(\gk,\gb_i,\nu)\in A^X_\nu$. Then by definition~\indiscDef(iic)
we have $c\in\cc^X(\gk,\gb)$ where $\gb>\gb_i$ for all sufficiently
large integers $i$.  If $c\not=s(\gk,\gb,\nu)$ then by
definition~\accptDef(2) it is because there is $\gb'\ge\gb$ and
$c'\in\cc(\gk,\gb')$ such that $\nu<c'<c$, but in this case we would
also have $c_i\not=s(\gk,\gb_i,\nu)$ for all $i$ large enough that
$c_i>c'$ and $\gb_i<\gb$.

\lemma  $A^{X}_{\nu } $
is cofinal in $\gk $ for all sufficiently large ordinals $\nu<\gk$.
\endth

\proof Suppose that the lemma is false and define a sequence $(\nu_k:k\in\gw)$
of ordinals less than $\gk$ as follows: set $\nu_0=0$ and for $k>0$
let $\nu_k$ be  the least member of
$X\setminus\sup(A^{X}_{\nu_{k-1}})$.   Let $\ga=\sup_k\nu_k$, so that
$\ga<\gk$ since $\cof(\gk)>\gw$. Let $\gb_k\le o^m(\gk)$ be the least
ordinal $\gb$ such that
$\union_{\gl\ge\gb}\cc(\gk,\gl)\setminus\nu_k=\nothing$.  Then
$\gb_{k+1}\le\gb_k$ so there are $\gb$ and $k_0$ such that $\gb_k=\gb$
for all $k\ge k_0$.  Then definition~\indiscDef ii(c) implies that
$\gb^{m}(\ga)\ge\gb$, contradicting  the definition of $\gb_{k_0}$.
\endproof

For each $X\prec H_{\gk^{++}}$ let $\nu^{X}$ be the least ordinal $\nu$ such
that $A^{X}_{\nu}$ is cofinal in $\gk $ and let $C^{X} =A^{X}_{\nu^X}$.  We
need a notion of supports which is more general than that of
definition~\accptDef.
Fix $\eta^X$ in the interval $\nu^X\le\eta^X<\gk$ so that the set
$$D^X=\set{c\in C^{X}:\cof(c)=\gw\text{ and }\gb^{X}(c)\in
h^{X}\image(\eta^X)}$$ 
is stationary. 
Call a finite sequence $d$ of ordinals an
{\it extended $X$-support} if for each $d_i\in d$ one of the following
two cases holds: 
\roster
\item"(i)" There are $\ga<\gk$ and  $\gb<o(\ga)$ in $h^X\image
(\eta^X\cup d\restrict i)$ 
such that either $d_{i}=s(\ga,\gb,d_{i-1})$, or $d_{i}=a(\ga,\gb,d_{i-1})$, or
\item"(ii)" For every $\ga\in D^X\setminus d_i+1$ there is $\gb\in
h^X\image(d\restrict i\cup\{\ga\}\cup\eta^X)$ such that either
$d_i=s(\ga,\gb,d_{i-1})$ or $d_i=a(\ga,\gb,d_{i-1})$.
\endroster

As in the last section, we say that $d$ is an extended $X$-support for $e$ if
$d$ is an extended $X$-support and every member $\nu$ of $e$ is either in $d$
or in $h^X\image(\eta^X\cup(d\cap\nu))$. 

\lemma For every finite sequence $e$ of ordinals in $X$ there is a sequence
$d$ which is an extended $X$-support for $e$.
\endth

\proof
Define $\gd_\ga(\gl)=\cohere^m(\gk,\gl,\gb(\ga))(\ga)$ for $\ga\in
D^X$ and $\gl<\gb(\ga)$. We will extend the inductive proof of
lemma~\havesupports\ by showing that for each
$c\in\union_{\gb}\cc(\gk,\gb)$ there are ordinals $\gl\in h\image
c\cap o^m(\gk)$ and $\gg<c$ such that either
$c=s(\ga,\gd_\ga(\gl),\gg)$ for all $\ga\in D^X\setminus c+1$ or
$c=a(\ga,\gd_\ga(\gl),\gg)$
for all $\ga\in D^X\setminus c+1$.
Since $\gl\in h\image c$ there is a finite subset $q$ of $c$ such
that $\gl\in h\image q$, and since $\gb(\ga)\in h\image\eta^X$ it
follows that $\gd_\ga(\gl)\in h\image(\eta^X\cup q\cup\{\ga\})$.  By
the induction hypothesis there is an extended $X$-support $d'$ for
$q\cup\{\gg\}$, and then $d=d'\cup\{c\}$ is an extended support for
$c$, completing the proof of the lemma.

The proof depends on the fact that $\gb(\ga)<\eta^X$ for all $\ga\in
D^X$, which implies that $\gd_\ga(\gl)$ gives an exact
correlation between indiscernibles for $\gk$ and those for ordinals
$\ga\in D^X$: if $\cc(\gk,\gl)\cap[\eta^X,\ga)\not=\nothing$ then
$\gd_\ga(\gl)$ is 
defined and $$\cc(\gk,\gl)\cap[\eta^X,\ga)=
\cc(\ga,\gd_\ga(\gl))\setminus\eta^X,\tag1$$ and conversely if
$\cc(\ga,\gd)\not=\nothing$ then there is $\gl$ such that
$\gd=\gd_\ga(\gl)$ and equation~(1) holds.

Now if $\ga\in D^X\setminus c+1$ then $c\in\cc(\ga,\gd_\ga(\gb(c)))$ and
it follows from lemma~\supportOne\ that for some $\gg<c$ either
$c=s(\ga,\gd_\ga(\gb(c)),\gg)$ or there is $\gd$ such that
$c=a(\ga,\gd,\gg)$. Suppose first that $c=s(\ga,\gd_\ga(\gb(c)),\gg)$.
Then for any $\ga'$ in $D^X\setminus c$ we have
$$
\union_{\gb\ge\gd_{\ga'}(\gb(c))}\cc(\ga',\gb)\cap(\eta^X,c]=
\union_{\gb\ge\gb(c)}\cc(\gk,\gb)\cap(\eta^X,c]
=\union_{\gb\ge\gd_{\ga}(\gb(c))}\cc(\ga,\gb)\cap(\eta^X,c]
$$ and it 
follows that $c=s(\ga',\gd_{\ga'}(\gb(c)),\gg)$ for all $\ga'\in
D^X\setminus c+1$, as required.

Thus we can assume that for each $\ga\in D^X\setminus c+1$ there is
$\gl_\ga\in h\image(c\cup\{\ga\})$ and $\gg_\ga$ such that
$c=a(\ga,\gl_\ga,\gg_\ga)$.  Since $\gl_\ga\in h\image(c\cup\{\ga\})$,
there is $\gl'_\ga$ such that $\gl_\ga=\gd_\ga(\gl'_\ga)$. Now
$c=a(\ga,\gl_\ga,\gg_\ga)$ implies that
$\union_{\gb\ge\gl_a}\cc(\ga,\gb)$ is bounded in $c$.  If $\ga'$ is
any other member of $D^X\setminus c+1$ then the argument in the last
paragraph implies that
$\union_{\gb\ge\gd_{\ga'}(\gl'_{\ga})}\cc(\ga',\gb)$ is bounded in
$c$, so that $c$ is not an accumlation point for $(\ga',\gb)$ for any
$\gb>\gd_{\ga'}(\gl'_\ga)$. Since $c$ is an accumulation point for
$\bigl(\ga',\gd_{\ga'}(\gl'_{\ga'})\bigr)$ it follows that
$\gl'_{\ga'}\le\gl'_{\ga}$, and since $\ga$ and $\ga'$ were arbitrary
it follows that $\gl'_\ga$  does not vary with $\ga\in D^X\setminus
c+1$. A similar argument shows that $\gg_\ga$ also does not vary with
$\ga\in D^X\setminus c+1$, and this completes the proof of the lemma.
\endproof

\demo{Proof of lemma~\mainthree}
 Fix $w\subset\gk $ of cardinality less than $\gk $, together with a
set $X\prec H_{\gk^{++}}$ such that  $w\subset X$, and let
$h=h^{X}\restrict\set{\nu:h^{X}(\nu)\in H_{\gk}}$.  Define the
increasing sequence $(\ga_{i}:i<\gw)$ of ordinals below $\gk$ by
induction on $i$: First set $\ga_{0}$ equal to $\eta^X$. Now suppose
that $\ga_{i}$ has been defined, and set $X_{i} =X^{h,\ga_i}$.  By
corollary~\uniqueFinally\ there is an ordinal $\xi_i<\gk$ such that
for any ordinals  $\ga,\gb,\gg\in X\cap X_{i} $ with $\gg>\xi_i$ we
have $s^{X}(\ga,\gb,\gg)=s^{X_i}(\ga,\gb,\gg)$ and
$a^{X}(\ga,\gb,\gg)=a^{X_i}(\ga,\gb,\gg)$.  Let $\ga_{i+1}$ be an
ordinal such that $\ga_{i+1}>\max\{\xi_i,\ga_i\}$, and $\ga_{i+1}$ is
sufficiently large that $\set{c\in D^X\cap C^{X_i}:\gb^{X_i}(c)\in
h^{X_i}\image(\ga_{i+1})}$ is stationary. We claim that
$X\subset\union_{i} X_{i} $.  By lemma~\backthm0 it is enough to show
that every extended $X$-support $d$ is in $\union_i X_i$. Let
$\nu=\sup_{i}\ga_{i}$ and let $i$ be large enough that
$d\cap\nu\subset\ga_{i}$.  We will show that $d\in X_{i}$.

Since $d\cap\nu\subset\ga_{i}\subset X_i $ it is enough to prove by
induction on $j$ that $d_{j}\in X_i$ for all $j$ with $d_j\ge\nu$.
Suppose that
$d\restrict j\in X_i$. Now if case~(i) of the definition of an extended
$X$-support holds for $d_j$ then $d_{j}$ is equal to either $s^{X}(\ga,\gb,d_{j-1})$
or $a^{X}(\ga,\gb,d_{j-1})$
for some $\ga$ and $\gb $ in $h^{X}\image (\eta^{X}\cup d\restrict j)$.
Since $\eta^X\cup d\restrict j$ is in $X_i$ we know that $\ga$ and $\gb$ are
in $X_i$, and
hence $s^{X_i}(\ga,\gb,d_{j-1})$ is in $X_i$.  
Now $s^X(\ga,\gb,\gg)=s^{X_i}(\ga,\gb, \gg)$ (and similarly for
$a(\ga,\gb,\gg)$) 
for all $\gg\ge\ga_i$, so we are done unless $d_{j-1}<\ga_i$ in which case
$d_j=s^X(\ga,\gb,d_{j-1})=s^X(\ga,\gb,\ga_i)=s^{X_i}(\ga,\gb,\ga_i)\in X_i$
(or similarly for $a(\ga,\gb,d_{j-1})$).

Now if, on the other hand, case~(ii) holds then $\set{\ga\in
D^X:\gb^{X_i}(\ga)\in h^{X_i}\image\nu}$ is stationary, and hence is
unbounded.  Thus we can pick $\ga$ in this set so that $\ga>d_i$ and
use the same argument as for case~(i).
\endproof

\heading\newsectno4 Countable Cofinality\endheading

  The case of cofinality $\gw$ is more delicate than the case of uncountable
cofinality.  With $\cof(\gk)>\gw$ every countable sequence $\ve\ga$ of
ordinals below $\gk$ was bounded in $\gk$.  This meant that there were only
$\gk$ many such sequences and hence the theorem could be proved by covering
every small subset of $\gk$ with a union
$\union_kX^{h,\ga_k}$ of basic covering sets.
Since this trick won't work when $\cof(\gk)=\gw$ we will have to define
standard systems of indiscernibles having a stronger maximality
property. These systems are given by lemma~\haveStdSys\ below.
In order to 
illustrate the method we begin with another result which uses a technique
suggested by Jensen to show that under a stronger hypothesis there is
system $\cc$ which does not depend at all on the set $X$ to be covered, and
hence is maximal in the strongest possible sense.  This generalizes the result
of Jensen and Dodd [D] that if $L(\mu)$ exists but $0^{\dagger} $ does not,
and there is a Prikry sequence over $L(\mu)$, then there is a unique (up to
finite changes) maximal Prikry sequence over $L(\mu)$. 

\theorem Assume that $\ga^{\gw }\le\ga^{+} $ for all ordinals $\ga$.  If there
is no regular limit of measurable cardinals in $\coref$ and no proper class
of measurable cardinals then there is a system $\cc$ of indiscernibles for
$\coref$ with the following maximality property: if $\cc'$ is any other
system of indiscernibles for $\coref$ then
$\union\set{\cc'(\xi)\setminus\cc(\xi):\xi\text{ is a cardinal}}$ is finite.
\endth 

\proof Since we are dealing with sequences $\cf$ which have at most one
measure per cardinal, we will write $\cf(\ga)$ instead of $\cf(\ga,0)$
and $\cc(\ga)$ instead of $\cc(\ga,0)$.

  We will use recursion over $\ga $ to construct systems $\cc_{\ga } $ of
indiscernibles for $\cf\restrict\ga $ which have the desired maximality
property.  Assume as an induction hypothesis that we have constructed
$\cc_{\ga'} $ for all $\ga'<\ga $.  The successor case is easy: a maximal
system of indiscernibles for $\cf\restrict\ga$ is also a maximal system for
$\cf\restrict(\ga+1)$ unless $\ga$ is measurable in $\coref$, in which case
the maximal system $\cc_{\ga +1} $ for $\cf\restrict(\ga+1)$ can be
obtained by adding a maximal Prikry sequence for $\cf(\ga)$ to the maximal
system $\cc_{\ga}$ for $\cf\restrict\ga$. 

  If $\ga$ is a limit ordinal but not a limit of measurable cardinals then
there is no problem, so we can assume that $\ga$ is a limit of measurable
cardinals.  Then the hypothesis implies that $\ga$ is singular in $\coref$,
say 
$\ga=\sup_{\nu <\gl }\ga_{\nu } $ where $\gl=\cof^{\coref}(\ga)$ and
$(\ga_{\nu }:\nu <\gl)\in \coref$ is continuous and unbounded in $\ga$. 
Take $X\prec H_{\ga^{++}}$ so that $\card X<\gk$, $\gl\subset X$, $\ga\in X$, and
$(\cc_{\xi}:\xi <\ga)\in X$.  The required system $\cc=\cc_{\ga}$ is
obtained by using $\ccx$ to combine the systems $\cc_{\ga_\nu}$:

$$\cc(\gz)=\cases\ccx(\gz) &\text{if $\gz\in X$}\\
\cc_{\ga_{\nu+1}}(\zeta) &\text{if $\zeta\notin X$ and
$\ga_{\nu } <\zeta<\ga_{\nu +1}$}.\endcases$$ 

  We claim that for each $\nu<\gl$ the restriction 
$\cc'=\cc\restrict (\ga_{\nu+1}\setminus\ga_{\nu})$ is a member of $X$.  
Consider the following two sets:
$$\alignat1
c_{1} &=\union
 \set{\ccx(\gz)\setminus\cc_{\ga_\nu}(\gz):\ga_{\nu}<\gz\le\ga_{\nu+1}}\\
c_2 &=\union
  \set{\cc_{\ga_\nu}(\zeta)\setminus\ccx(\zeta):\zeta\in X\text{ and }\ga_{\nu } <\zeta
	<\ga_{\nu +1} }
\endalignat$$
The set  $c_{1} $ is finite by the maximality of $\cc_{\ga_\nu}$, and since
$c_{1}\subset X$ it follows that $c_{1}\in X$.  
The set $c_2$ is finite by lemma~\unique. It is a subset of $X$: otherwise let
$\nu$ be the largest member  of $c_{2}\setminus X$ and and let $\gz$ be the
ordinal such that
$\nu\in\cc_{\ga_\nu}(\gz)$. Then $\cc_{\ga_{\nu}}(\gz)\in X$, and if $\gg$ is
the least member of $X$ above $\nu$ then $\nu$ is the largest member of
$\cc_{\ga_{\nu}}(\gz)$ below $\nu$ and hence is in $X$.
It follows that $c_2\in X$.
Then $\cc_{\ga_\nu}$, $c_{1} $ and $c_{2}$ are all in $X$ and since
$\cc_{\ga_{\nu}}$ can be converted to $\cc'$ by adding
$c_1$ and deleting $c_2$ it follows that  $\cc'\in X$. 

  Now suppose that $\cc$ is not a system of indiscernibles.  Then there is
a countable set $I\subset\gl$ such that $\cc^I=\cc\restrict\union_{\nu\in I}
(\ga_{\nu +1}\setminus\ga_{\nu})$
is not a system of indiscernibles.  Since
${^\gw } X\subset X$, $\cc^I$ is in $X$ and so by elementarity it
is true in $X$ that $\cc^I$ is not a system of indiscernibles.  This
is absurd because $\cc\restrict X=\ccx$.  The same argument shows that
$\cc$ is maximal and completes the proof of the theorem. 
\endproof 

  A result due to Jensen and independently to myself and P.~Matet (see [Mi84a])
shows that the hypothesis of theorem~\backthm0\ 
cannot be weakened further: If there is a regular
limit $\gk $ of measurable cardinals or a class of measurable cardinals
then there is a model with the same cardinals in which every measurable
cardinal has a Prikry sequence, but in 
which there is no system $\cc$ of indiscernibles such that
$\cc(\ga)\not=\nothing $ for unboundedly many $\ga <\gk $.  

  Our proof of theorem~\mainthm(i) is based on the same ideas. We
will assume that the conclusion of theorem~\mainthm(i) is false, that
is, that there is an ordinal $\gb_0<\gk$ such that $o(\ga)<\gb_0$ for
all $\ga$, and use this to prove the SCH. Although there  is no
single maximal system of indiscernibles for $\cf$ we will be able to
construct a set of standard systems, each of which is a maximal
system of indiscernibles on domains determined by some particular
function in $\coref$.  We no longer are restricted to one measure per
cardinal, but since $o(\ga)<\ga$ for all $\ga$ some simplification of
notation is still possible. We write $\cc(\ga)$ for
$\union_{\gb<o(\ga)}\cc(\ga,\gb)$. Note that
$\cc(\ga,\gb)=\set{\nu\in\cc(\ga):o(\nu)=\gb}$.

  The restriction to cardinals in theorem~\backthm0 can be weakened slightly,
but it cannot be eliminated.  Thus we need to introduce some technical
apparatus so 
that we can restrict our attention to systems $\cc$ of indiscernibles 
such that $\cc(\ga)=\nothing$ whenever $\ga$ is not a cardinal in the real
world. 
For simplicity we will further restrict
ourselves to systems of indiscernibles such that $\cc(\ga)=\nothing$
whenever there is $\gl<\ga$ such that $\gl^\gw\ge\ga$.
Call an ordinal $\ga$  {\it full} if $\ga$ is a cardinal and $\gl^\gw<\ga$ for
all $\gl<\ga$.  
In order to
cover a set with such restricted systems fix functions $\gs $ and $\tau $
such that $\gs(\nu,\cdot):\card\nu\cong\nu $ and 
$\tau (\nu,\cdot):{^\gw}\nu\cong\card\nu^{\gw }$ for all ordinals $\nu $.   
For
the rest of this section we will assume that all systems of indiscernibles
mentioned are 
empty except on full cardinals and that every elementary substructure $X$ of
$H_{\gk^{++}}$ which we use contains the set $\{\gs,\tau\}\cup\gb_{0}$ and
also contains all of its limit points of cofinality at most $\gb_{0}$.
Notice that by corollary~\semicont\ this implies that there are at
most finitely many accumulation points in $\cc^X$, so  that we can
assume \iwlog\ that there are none.

 Let $\ch$ be the class of $h\in\coref$
such that there is a ordinal $\eta$ such that $h$
maps a cofinal subset of $\eta$ into $\eta$ so that $h(\nu)>\nu$ for all
$\nu\in\domain(h)$.  For $h\in\ch$ we will abuse notation by writing
$\domain(h)$ for $\eta$ and, if $\xi<\eta$, by writing $h\restrict\xi$ for
$h\cap\xi^2\in\ch$. 
The following definition is a modification of
notation from earlier in the paper: 

\defin\thmtag\hCclosureDef If $h\in\ch$ with $\eta=\domain(h)$ and
$\cc$ is a system of indiscernibles then we write
$h\image(\gg;\cc)$ for the smallest set $Y\subset\eta+1$ such that

(i)  $\gg\cup\{\eta\}\subset Y$ and $Y$ is closed under the canonical
pairing function,

(ii) $h\image Y\subset Y$, $\gs \image Y^{2}\subset Y$, and $\tau 
\image(Y\times [Y]^{\gw })\subset Y$, and 

(iii) $s^{\cc}(\ga,\gb,\gg')\in Y$ whenever $\ga,\gb,\gg'\in Y$ and
$s^{\cc}(\ga,\gb,\gg')<\eta$. \endth

Note that under this definition it is still true that if
$h=h^{X}\restrict\set{\nu:\nu<h^X(\nu)<\gk}$ then $X\cap\gk\subset
h\image(\rho^X;\ccx)$.

The  next  lemma  asserts  the  existence  of  the  standard  systems
of indiscernibles  which  are  used  in  the  proof  of
theorem~\mainthm(i).   We  will  first prove  theorem~\mainthm(i)
under  the  assumption  that  lemma~\haveStdSys\   holds,  and  then
prove  lemma~\haveStdSys. 

\lemma\thmtag\haveStdSys 
For all $h\in\ch$ there is a system $\cc^h$
such that for any
system $\cc$ of indiscernibles of the form $\cc^Y$, any sufficiently 
large $\gd<\eta=\domain(h)$, and any $\ga\in h\image(\gd;\cc^h)$
we have $\cc(\ga)\setminus\gd\subset\cc^{h}(\ga)$.  
\endth

The following lemma will be used in the proof of both theorem~\mainthm(i) and
of lemma~\haveStdSys.

\lemma\thmtag\XisMax Suppose $X\prec H_{\gk^{++}}$, $\cc$ is a system
of indiscernibles of the form $\cc^Y$,
$\gg\in X$, and $\gg >\rho^X$.  Then there is a $\gd<\gg$ such
that $\cc(\ga)\setminus\gd\subset\ccx(\ga)$ 
for all $\ga\in X\cap (\gg+1\setminus\gd)$.
\endth

\proof If either of $X$ or $\union_{\ga }\cc(\ga)$ is bounded in $\gg $
then the lemma is immediate.  If both are unbounded then by
elementarity $X$ satisfies that there is a system $\cd$ of
indiscernibles which is unbounded in $\gg$, and all but finitely many
of the indisceribles in $\cd$ must be indiscernibles in $\ccx$ since
$c\in h^X\image c$ for every ordinal $c$ which is not an
indiscernible in $\cc^X$. Thus $\union_{\ga}\ccx(\ga)$ is unbounded
in $\gg$.

Now suppose that $c$ and $\ga$ are ordinals such that
$c\in\cc(\ga)\setminus\cc^X(\ga)$. 
Set $\gd(c)=0$ if $\ccx(\ga)\cap c=\nothing$ and $\gd(c)=\sup(\ccx(\ga)\cap
c)$ otherwise, and set $\nu(c)=\inf(X\setminus c+1)$.  Then $\gd(c)$ and
$\nu(c)$
are both in $X$  (for $\gd(c)$, note that $\cof(\gd(c))\le\gb_0$ by
corollary~\semicont\ and recall that $X$ contains its limit points of
cofinality at most $\gb_0$) and
$\gd(c)<c<\nu(c)<\gg$ and $\ccx(\ga)\cap(\nu(c)\setminus\gd(c))=\nothing$.

Now if the lemma is false then there must exist
sequences $\vc$, $\ve\ga$, $\ve\gd$ and
$\ve\nu$ in $X$ such that for each $n\in\gw$ we have
$c_n\in\cc(\ga_n)\setminus\ccx(\ga)$, $\gd_n=\gd(c_n)$, $\nu_n=\nu(c_n)$
and $\nu_n\le\gd_{n+1}$.
Then
by elementarity $X$ satisfies that there is a sequence $\vd=(d_{n}:n\in \gw)$, a
system $\cd$ of indiscernibles, and a function $h\in \coref$ such that
for each $n$ we have $\gd_{n} <d_{n} <\nu_n$, $d_{n}\in\cc(\ga_{n})$, and
$\ga_{n}\in h\image d_{n} $.  Then $d_{n}\in\union_{\ga }\ccx(\ga)$ for all
but finitely many $n\in\gw $, for otherwise $d_{n}\in h^{X}\image(d_{n})$
for infinitely many $n$, contradicting the assumption that $\cd$ is a
system of indiscernibles.  Then $\ga_{n}\in h\image d_{n} $ implies that
$\ga^{X}(d_{n})\le\ga_{n} $ for all but finitely many $n$.  If
$\ga^{X}(d_{n})<\ga_{n} $ infinitely often then again $h^{X} $ would be a
counterexample to the indiscernibility of $d_{n}\in\cd(\ga_{n})$ in $V$,
and hence there would be such a counterexample in $X$.  Thus for infinitely
many $n$ we have $\ga_{n}=\ga^{X}(d_{n})$ and hence $d_n\in\ccx(\ga_{n})$,
contradicting the definition of $\gd_{n} $ and $\nu_n$.  
\endproof

\demo{Proof of theorem~\mainthm({\rm i}) assuming lemma~\haveStdSys} 
Suppose
that $X\prec H_{\gk^{++}}$ is a set as used in the covering lemma.  Let
$m=m^X$ be the covering mouse, let $h=h^{m^X}\cap\gk^2$,
and let $\rho$ be large enough that $X\cap\gk\subset h\image(\rho;\cc^X)$.  
We will assume that corollary~\semicont\ holds for every ordinal $c$ in 
$\cc^X$, increasing $\rho$ if necessary to be larger
than any of the finitely many exceptions.

We will find an ordinal $\gd<\gk$ such that $X\subset
Y=Y^{\gd,h}=h\image(\gd;\cc^h)$.  Since there 
are only 
$\gk^+$ functions $h:\gk@>>>\gk$ in $\coref$ and only $\gk$ ordinals
$\gd<\gk$ there are only $\gk^+$ sets $Y^{\gd,h}$.
Since each set $Y^{\gd,h}$ has cardinality
$\card\gd<\gk$ there are fewer than $\gk$ subsets of $Y^{\gd,h}$.  
Thus there are only $\gk^+$ subsets of $\gk$ of cardinality less then $\gk$
and hence $2^{\gk}=\gk^+$.

Let $\gd$ be any ordinal in the interval $\rho<\gd<\gk$ such that
$\ccx(\ga)\setminus\gd=\cc^h(\ga)\setminus\gd$ for all $\ga\in
h\image(\gd;\cc^h)\cap X$.  We can always find such a $\gd$ since
lemma~\haveStdSys\ asserts that for sufficiently large $\gd<\gk$ we have
$\ccx(\ga)\setminus\gd\subset\cc^h(\ga)$ for all $\ga\in h\image(\gd;\cc^h)$,
while lemma~\XisMax\ asserts that for sufficiently large $\gd<\gk$ we have
$\ccx(\ga)\setminus\gd\supset\cc^h(\ga)$ for all $\ga\in X$. 

Now suppose $\nu\in(\gk\cap X)$.  We will show by induction on $\nu$ that
$\nu\in Y=h\image(\gd;\cc^h)$.  If $\nu<\gd$ then certainly $\nu\in Y$.  If
$\nu$ is not a full cardinal then either $\nu\in\gs\image ((X\cap\nu)^{2})$ or
$\nu\in\tau\image ((X\cap\nu)\times [X\cap\nu ]^{\gw })$ and since $X\cap
\nu\subset Y$ by the induction hypothesis $\nu$ is in $Y$ by
definition~\hCclosureDef(ii) above.  If $\nu$ is a full cardinal, but not an
indiscernible in $\ccx$, then $\nu\in h^{X} \image(\nu\cap
X)\cap\gk=h\image(\nu\cap X)\subset Y$, since $Y$ is closed under $h$.  Thus
we can assume that $\nu$ is an indiscernible in $\ccx$, so that
$\nu=s^{\cc^X}(\ga,\gb,\gg)$ for some $\gg\in X\cap\nu $, some $\ga\in
h^{X}\image(\nu\cap X)$, and some $\gb <\gb_{0} $.  Then $\gg\in Y$ by the
induction hypothesis.  The set $Y$ is closed under $h^{X} $ since it is closed
under $h$.  Hence $\ga $ and $\gb $ are in $Y$ and by the choice of $\gd$ we
have $\nu=s^{\cc^X}(\ga,\gb,\gg)=s^{\cc^h}(\ga,\gb,\gg)\in Y$.  \endproof

This completes the proof of theorem~\mainthm(i)\ from  lemma~\haveStdSys.
The proof of lemma~\haveStdSys\ will be by induction on the domain $\eta$ of
the function $h$.  The next proposition will be needed for the induction
step. 

\prop
Suppose that $h\in\ch$ and 
$\cc^f$ has been already defined as required by lemma~\haveStdSys\
for all functions $f\in\ch$ with
$\domain(f)\le\eta=\domain(h)$.  Assume further that 
$X\prec H_{\gk^{++}}$, the sequence 
$\seq{\cc^{f}:f\in\ch\mathrel{\text{and}}\domain(f)\le\eta}$ is a member of
$X$, and $\rho<\eta$.  
Then there is a system $\cc\in X$ such that
$\cc\restrict X=\ccx\restrict (h\image(\rho;\ccx)\cap((\eta+1)\setminus\rho))$.
\endth

\proof The proof is by induction on $\eta$.  If $h\image(\rho;\ccx)$ is
bounded in $\eta $ then the proposition follows immediately from the induction
hypothesis, so we can assume that $h\image(\rho;\ccx)$ is unbounded in
$\eta$. 

  Now lemma~\XisMax\ implies that there is $\gd_0<\eta$ such that 
$\ccx(\gg)\supset\cc^h(\gg)\setminus\gd_0$ for all $\gg\in X$, and the
hypothesis to this proposition implies that there is $\gd_1$ such that
$\ccx(\gg)\setminus\gd_1\subset\cc^h(\gg)$ for all 
$\gg\in h\image(\gz;\ccx)$.  
Let $\gd=\max(\gd_0,\gd_1)$.

Let
$\vc=
\set{c\le\gd:\exists\ga>\gd\;c\in\cc^X(\ga)
\mathrel{\text{and}}\forall\ga\le\gd\;c\notin\ccx(\ga)}$.  Then $\vc$ is
finite by definition~\indiscDef ii(c).  Define 
$h^{\vc}\in\ch$ with $\domain(h)=\gd$ by setting $h^{\vc}(i)=c_i$ 
for $i<\len(\vc)$ and $h^{\vc}(\len(\vc)+\nu)=h(\nu)$ otherwise.  
The induction hypothesis implies that there is a system $\cc_0\in X$ of
indiscernibles such that
$\cc_0\restrict X=
\ccx\restrict\big(h^{\vc}\image(\rho,\ccx)\cap(\gd+1\setminus\rho)\big)$.
Now define $\cc_1$ by using $\vc=(c_i:i<k)$ to combine $\cc_0$ and
$\cc^{h}\restrict (\eta\setminus\gd)$. That is, 

$$
\cc_1(\nu)=\cases\cc_0(\nu)&\text{if $\rho<\nu\le\gd$}\\
\cc^h(\nu)\cup\union\set{\cc_0(c_{i})\cup\{c_{i}\}:c_{i}\in\ccx(\nu)} 
&\text{if $\gd<\nu\le\eta$.}\endcases
$$

The system $\cc=\cc_1\restrict h\image(\rho ;\cc_1)$ is the required system.
\endproof

\demo{Proof of lemma~\haveStdSys} We will construct the systems $\cc^h$ by
induction over $\eta=\domain(h)$.  We assume as an induction hypothesis that
$\cc^f$ has been defined for all $f$ with $\domain(f)<\eta$.  If $\eta$ is
regular then $h\image(\rho;\cc)$ cannot be cofinal in $\eta$, so we can define
$\cc^h(\nu)=\nothing$ for all $\nu$.  Now assume that $\eta$ is singular.  If
$\eta$ is singular in $\coref$ then set $\gl=\cof^{\coref}(\eta)$, and
otherwise let $\gl=\cof(\eta)$.  
Take $X\prec H_{\eta^{++}}$ so that $X$ is cofinal in $\eta$, the sequence
$\seq{\cc^f:\domain(f)<\ga}$ is in $X$, and $\gb_0\cup\{p,\gb_0,\gl\}\subset X$.
If $\eta$ is singular in $\coref$ then let $(\gg_\gi:\gi<\gl)\in
\coref\cap X$ be a closed unbounded subset of $\eta$. 
Otherwise $\ccx(\eta)$ is unbounded in $\eta$, so let $(\gg_\gi:\gi<\gl)$ be a
$\gw$-closed, unbounded subset of $\ccx(\eta)$.  For each $\gi<\gl$ let
$h_\gi=h\restrict\gg_{\gi+1}$ and use proposition~\backthm0 to pick
$\cc^{\gi}\in X$ so that $\cc^{\gi}\restrict X=\ccx\restrict
\big(h_\gi(\gg_{\gi};\ccx)\cap(\gg_{\gi+1}+1\setminus\gg_{\gi})\big)$.  Define
$\cc^h(\nu)=\cc^{\gi}(\nu)$ for $\gg_{\gi}<\nu\le\gg_{\gi+1}$, and if $\eta$
is measurable in $\coref$ then set $\cc^h(\eta)=
\set{\gg_{\gi}:\gi<\gl}\cup\union_{\gi}\cc^{\gi}(\gg_{\gi+1})$. 

  The proof that $\cc^h$ is a system of indiscernibles is just like that
in theorem~4.1: if it is not, then there is a countable set
$I\subset\gl$ such that  $\cc^I=\cc^h\restrict
\union_{\gi\in I}(\gg_{\gi +1}+1\setminus\gg_\gi)$ is not
a system of indiscernibles.  But $\cc^I$  is a member of
$X$ so it is true in $X$ that $\cc^I$ is not a system of indiscernibles, which
is impossible  since $\cc^I\restrict X$ is the restriction of $\ccx$ to the
set $\union_{\gi\in I}(\gg_{\gi+1}+1\setminus\gg_\gi)$.  Similarly if
$\cc^h$ does not satisify the maximality property of lemma~\haveStdSys\
then there
must be another  system $\cc$ of indiscernibles
and a countable set $I\subset\gl$ such that for
each $\gi\in I$ there is
$\nu_{\gi}\in h\image(\gg_{\gi};\cc)\cap(\gg_i,\gg_{\gi+1}+1]$ such
that $\cc(\nu_{\gi})\setminus\gg_{\gi}$ is not contained in
$\cc^I(\nu_{\gi})$.  But then there must be such a system $\cc$ and
sequence $\nu$
which is a member of $X$, so that except for finitely many $\gi$ we must have
$\cc(\nu_{\gi})\subset\ccx(\nu_{\gi })=\cc^h(\nu_{\gi })$, a contradiction. 
\endproof

\medskip
I would like to thank the referee for a careful reading of this paper
and many helpful suggestions.

%________________________________________________________________

\enddocument